\documentclass[twocolumn,10pt]{IEEEtran} 
\usepackage{amsfonts}
\usepackage{amssymb}
\usepackage{mathrsfs} 
\usepackage{amsmath}
\usepackage{amsfonts, amssymb}
\newtheorem{theorem}{Theorem}
\usepackage{graphicx}
\usepackage{color}
\usepackage{algorithm}
\usepackage{algorithmic}

\newtheorem{corollary}{Corollary}
\newtheorem{definition}{Definition}
\newtheorem{lemma}{Lemma}
\newtheorem{example}{Example}

\newtheorem{problem}{Problem}
\newtheorem{assumption}{Assumption}
\newtheorem{remark}{Remark}
\renewcommand{\baselinestretch}{1.19}
\def\BibTeX{{\rm B\kern-.05em{\sc i\kern-.025em b}\kern-.08em
		T\kern-.1667em\lower.7ex\hbox{E}\kern-.125emX}}
\allowdisplaybreaks
\DeclareMathOperator{\rank}{rank}
\DeclareMathOperator{\Null}{Null}

\usepackage[numbers, sort&compress]{natbib}

\makeatletter
\newcommand{\rmnum}[1]{\romannumeral #1}
\newcommand{\Rmnum}[1]{\expandafter\@slowromancap\romannumeral #1@}
\makeatother

\begin{document}
	
	\title{Signed Angle Rigid Graphs for Network Localization and Formation Control}
	\author{Jinpeng Huang~and~Gangshan Jing
    \thanks{J.~Huang and G.~Jing are with School of Automation, Chongqing University, Chongqing 400044, China,
		{\tt\small 202313021067t@stu.cqu.edu.cn, jinggangshan@cqu.edu.cn}}}
    
	
	\markboth{{\it draft}} {Hu: Using the style file IEEEtran.sty} \maketitle
	
	\begin{abstract}
		Graph rigidity theory studies the capability of a graph embedded in the Euclidean space to constrain its global geometric shape via local constraints among nodes and edges, and has been widely exploited in network localization and formation control. In recent years, the traditional rigidity theory has been extended by considering new types of local constraints such as bearing, angle, ratio of distance, etc. Among them, the signed angle constraint has received extensive attention, since it is practically measurable and independent of the global coordinate frame. However, the relevant studies always consider special graph structures, which are sufficient but not necessary for signed angle rigidity. This paper presents a comprehensive combinatorial analysis in terms of graphs and angle index sets for signed angle rigidity. We show that Laman graphs equivalently characterize minimally signed angle rigid graphs. Moreover, we propose a method to construct the minimal set of signed angle constraints in a Laman graph to effectively ensure signed angle rigidity. These results are finally applied to distributed network localization and formation stabilization problems, respectively, where each agent only has access to signed angle measurements.
        
	\end{abstract}
	\section{Introduction}\label{sec:1}

    The traditional graph rigidity theory studies the characterizability of a geometric shape by distance constraints among neighboring vertices in the graph, and has been widely employed in multi-agent coordination problems such as network localization \cite{aspnes2006theory,fang2009sequential, wan2019sensor} and formation control \cite{anderson2008rigid, krick2009stabilisation, lin2015necessary,zhao2018affine,zhang2025linear}. Inspired by this, rigidity theories based on other types of constraints (such as bearings \cite{eren2012formation, zhao2015bearing}, angles \cite{jing2019angle, chen2020angle}, ratio of distance \cite{cao2019ratio}, etc.) have been further developed to achieve network localization \cite{fang2020angle,cano2023ranging,zhao2016localizability,jing2021angle,chen2022triangular} and formation stabilization \cite{eren2012formation, zhao2015bearing, jing2019angle,chen2020angle, wang2024adaptive,cao2023similar, cao2019ratio} via different sensing measurements. A geometric shape specified by different types of constraints may require different graphical conditions to ensure a unique shape \cite{oh2015survey}. As a result, different conditions on the network topology have been proposed depending on the type of measurements utilized.

    In distance-based multi-agent coordination problems, infinitesimal distance rigidity has been adopted most widely since it has an equivalent rank condition on the rigidity matrix \cite{asimow1978rigidity}, which plays an essential role in the analysis of the gradient-based controller \cite{olfati2002graph,krick2009stabilisation,de2016distributed, chen2023distance}. However, infinitesimal distance rigidity does not guarantee the shape to be uniquely determined by distance constraints \cite{jackson2007notes}, which usually leads to local convergence of the controller. On the other hand, although global distance rigidity implies a unique shape, and has been favored in sensor network localization \cite{aspnes2006theory}, it usually requires a large number of edges. Moreover, the formation controller achieving global convergence under globally rigid graphs is yet to be designed \cite{park2017distance}. In \cite{zhao2015bearing}, infinitesimal bearing rigidity was shown to share the same graphical conditions with infinitesimal distance rigidity (i.e., Laman graphs \cite{zhao2017laman,jackson2007notes}), while guaranteeing that the geometric shape is uniquely determined by bearings. Guided by this observation, the authors in \cite{zhao2015bearing, zhao2016localizability} proposed novel distributed protocols for bearing-based formation stabilization and sensor network localization. However, bearing constraints are associated with the global coordinate frame, which limits its applicability.

    Recently, angle rigidity theory has garnered significant interest, as angle measurements are independent of the coordinate frame \cite{fang2020angle, peng2025angle, buckley2021infinitesimal, jing2021angle, chen2022triangular, jing2018weak}, and an angle-constrained formation exhibits higher degrees of freedom than distance and bearing-based approaches \cite{jing2019angle,chen2020angle}. In \cite{jing2019angle}, the authors established the equivalence between infinitesimal angle rigidity, infinitesimal bearing rigidity, and infinitesimal distance rigidity in the plane. However, it is unclear whether the shape of an infinitesimally angle rigid framework can be uniquely determined by angles. In \cite{chen2020angle}, the authors reformulated angle rigidity theory by considering signed angle constraints. Unfortunately, the signed angle rigidity property was studied on angularities with predefined sets of angle constraints. Such a setting benefits for excluding redundant angles, but in the meanwhile, hinders the connection between the graph structure and the signed angle rigidity property. 

    In this paper, we study the conditions on graphs and the minimal set of signed angle constraints required to ensure signed angle rigidity and uniqueness of the shape. Similar issues related to angle rigidity theory have been investigated in the literature. In \cite{zhou2006combinatorial,dewar2024angular}, the authors explored the combinatorial conditions of angle rigidity. In \cite{chen2022triangular, chen2020angle, jing2019angle}, the authors reported a class of particular graphs which ensure a unique shape under angle constraints. Nevertheless, these works fail to characterize the equivalent graphical conditions for angle rigidity and the minimal set of angle constraints to guarantee a unique shape. Considering the large body of existing works in bearing rigidity theory, we believe that establishing the relationship between signed angle rigidity and bearing rigidity would be a worthwhile exploration. 

    Our contributions can be summarized as follows. (\rmnum{1}). Enlightened by angle rigidity theory in \cite{jing2019angle, chen2020angle}, we develop the signed angle rigidity theory in the context of frameworks, where all signed angles within the framework are utilized. Under this setting, it is shown that (infinitesimal, global) signed angle rigidity and (infinitesimal, global) bearing rigidity are completely equivalent, which implies that Laman graphs are exactly minimally signed angle rigid graphs; see Section \ref{sec:2}--\ref{sec:3}. (\rmnum{2}). Similarly to the bearing case, we prove that infinitesimal signed angle rigidity is sufficient and necessary for a framework to determine its shape uniquely; see Theorem \ref{theorem shape fix}. (\rmnum{3}). We provide a combinatorial analysis on the set of angle constraints and develop Algorithm \ref{alg:1} to generate a minimal angle index set (AIS) to characterize the rigidity property of a framework; see Section \ref{sec:4}. (\rmnum{4}). Based on the theory developed, we solve the problems of network localization and formation stabilization, respectively, where each agent only has access to local signed angle measurements associated with its neighbors.  Two distributed algorithms with global convergence guarantee are proposed, respectively, where the conditions on the sensing graph are much milder than those in \cite{chen2022triangular,jing2019angle,jing2021angle,chen2020angle}; see Section \ref{sec:5}--\ref{sec:6}.

    This paper is structured as follows. Section \ref{sec:2} reviews bearing rigidity theory and introduces the concept of signed angle rigidity theory for frameworks in the plane. Section \ref{sec:3} establishes the equivalence on bearing rigidity and signed angle rigidity. Section \ref{sec:4} provides a combinatorial analysis for minimal angle index sets. Section \ref{sec:5}-\ref{sec:6} demonstrate the application of our results to network localization and formation stabilization. Section \ref{sec:7} summarizes the entire paper.
    
	Notations: In this paper, we denote the set of real numbers by ${\mathbb{R}}$; Let ${\mathbb{R}}^{d}$ be the $d$--dimensional Euclidean space; $\|\cdot\|$ is the Euclidean norm; $|A|$ indicates the number of the elements in the set $A$; $X^{\top}$ denotes the transpose of the matrix $X$; $\otimes$ is the Kronecker product; Null($X$) and rank($X$) are defined as the null space and rank of the matrix $X$, respectively; Define $I_n \in \mathbb{R}^{n\times n}$ as the identity matrix; $\textbf{1}_n \triangleq (1,\cdots,1)^{\top}$; Denote $A\backslash B$ as the set that consists of the elements belonging to $A$ but not to $B$; Denote  $\mathscr{R}_o(\theta)$ as the 2-dimensional rotation matrix associated with $\theta \in [0,2\pi)$; For a vector $x\in {\mathbb{R}}^{2}$, $x^\perp \triangleq \mathscr{R}_o(\frac{\pi}{2})x$; For $X_{i}\in {\mathbb{R}}^{a\times b}, i = 1,\cdots,q$, diag$\left(X_i\right) \triangleq $ blockdiag$\left\{X_1,\cdots,X_q\right\}\in{\mathbb{R}}^{qa\times qb}$.
	
    We consider the undirected graph ${\mathcal{G}} \triangleq \left(\mathcal{V}, \mathcal{E}\right)$, where $\mathcal{V} \triangleq \left\{1,\cdots,n\right\}$ is the vertex set satisfying $n\geq3$, $\mathcal{E} \subseteq \mathcal{V} \times \mathcal{V}$ is the edge set, and $(i,j)\in\mathcal{E}$ is equivalent to $(j,i)\in\mathcal{E}$ in $\mathcal{G}$. The set of neighbors of the vertex $j$ in $\mathcal{G}$ is denoted by $\mathcal{N}^j_{\mathcal{G}} \triangleq \left\{j\in\mathcal{V}:(j,i)\in\mathcal{E}\right\}$. If $\mathcal{G}$ is clear from the context, we simply denote $\mathcal{N}^j_{\mathcal{G}}$ by $\mathcal{N}_j$. The undirected complete graph with the same number of nodes as $\mathcal{G}$ is denoted by $\mathcal{K}$. A graph $\bar{\mathcal{G}} \triangleq (\bar{\mathcal{V}},\bar{\mathcal{E}})$ is said to be a subgraph of $\mathcal{G}$ if $\bar{\mathcal{V}} \subseteq \mathcal{V}$ and $\bar{\mathcal{E}} \subseteq (\bar{\mathcal{V}} \times \bar{\mathcal{V}}) \bigcap \mathcal{E}$. A subgraph is said to be a spanning subgraph if it is connected and $\bar{\mathcal{V}} = \mathcal{V}$. 

    Given an undirected graph $\mathcal{G}$, an undirected path between $l_1\in\mathcal{V}$ and $l_{n_P}\in\mathcal{V}$ is a sequence of edges of the form $(l_1,l_2)$,$\dots, (l_{n_P-1},l_{n_P})$. In this paper, for the sake of analysis, we use the subgraph $\mathcal{P} \triangleq (\mathcal{V}_P,\mathcal{E}_P)$ to represent the corresponding undirected path, where $\mathcal{V}_P \triangleq \{l_1,\dots,l_{n_P}\} \subseteq \mathcal{V}$ and $\mathcal{E}_P \triangleq \{(l_1,l_2),\dots,(l_{n_P-1},l_{n_P})\} \subseteq\mathcal{E}$. An undirected cycle is an undirected path in $\mathcal{G}$ that starts and ends at the same node. We use the subgraph $\mathcal{C} \triangleq(\mathcal{V}_C,\mathcal{E}_C)$ to represent the corresponding undirected cycle, where $\mathcal{V}_C \triangleq \{l_1,\dots,l_{n_P}\} \subseteq \mathcal{V}$ and $\mathcal{E}_C \triangleq \{(l_1,l_2),\dots,(l_{n_P-1},l_{n_P}),(l_{n_P},l_1)\} \subseteq\mathcal{E}$.
    
    \section{Preliminaries}\label{sec:2}
	
	In this section, we briefly review some background on bearing rigidity theory \cite{zhao2015bearing}, define signed angle rigidity in the context of frameworks, and introduce some concepts related to the index set of angle constraints. 
	
	\subsection{Graph Rigidity Theory}
	
	To begin with, we introduce some important notions in rigidity theory based on a smooth rigidity function. Given different smooth rigidity functions, different rigidity theories can be developed accordingly, e.g., bearing rigidity theory and signed angle rigidity theory.
	
	A framework $(\mathcal{G}, p)$ in $\mathbb{R}^d$ is characterized by an undirected graph $\mathcal{G}$ and a configuration $p = \left(p^{\top}_1, \cdots, p^{\top}_n\right)^{\top}\in\mathbb{R}^{dn}$, where $p_i$ is the coordinate mapped to the vertex $i, i = 1,\cdots, n$ with the assumption that $p_i \neq p_j$ for any $j\in\{1,\cdots, n\}\backslash\{i\}$. If $p_1,\cdots, p_n$ do not belong to a hyperplane of $\mathbb{R}^{d}$, we refer $p$ to be \textit{non-degenerate}, and the framework $(\mathcal{G}, p)$ is said to be a \textit{non-degenerate framework}. Given a spanning subgraph $\bar{\mathcal{G}}$ of $\mathcal{G}$, $(\bar{\mathcal{G}}, p)$ is said to be the \textit{spanning subframework} of $(\mathcal{G},p)$.
	
	Given a framework $(\mathcal{G}, p)$, a smooth \textit{rigidity function} $r_{\mathcal{G}}(\cdot): \mathbb{R}^{dn} \rightarrow \mathbb{R}^{s}$ is a function of $p\in\mathbb{R}^{dn}$, and each entry of $r_\mathcal{G}(p)$ corresponds to a constraint in $(\mathcal{G}, p)$, here $s$ is the number of constraints. Based on $r_{\mathcal{G}}(\cdot)$, a series of relevant concepts can be given.
	
	\begin{definition}
		In $\mathbb{R}^d$, a framework $(\mathcal{G}, p)$ is said to be \textit{rigid} if there exists a neighborhood $U_p$ of $p$, such that for any $q \in r^{-1}_{\mathcal{G}}\left(r_{\mathcal{G}}(p)\right)\bigcap U_p$, it always holds $r_{\mathcal{K}}(q) = r_{\mathcal{K}}(p)$. 
	\end{definition}
	
	\begin{definition}
		In $\mathbb{R}^d$, a framework $(\mathcal{G}, p)$ is said to be \textit{globally rigid} if $r^{-1}_{\mathcal{G}}\left(r_{\mathcal{G}}(p)\right) = r^{-1}_{\mathcal{K}}\left(r_{\mathcal{K}}(p)\right)$.
	\end{definition}

        \begin{definition}
            In $\mathbb{R}^d$, a rigid framework $(\mathcal{G}, p)$ is said to be $\textit{minimally rigid}$ if there does not exist any rigid spanning subframework $(\bar{\mathcal{G}},p)$ of $(\mathcal{G}, p)$ such that $|\bar{\mathcal{E}}|<|\mathcal{E}|$.
        \end{definition}
	
	By definition, global rigidity implies rigidity.
	
	An instantaneously continuous motion $\dot{p}\triangleq v = (v^{\top}_1, \cdots, v^{\top}_n)^{\top}$ is said to be an \textit{infinitesimal motion} if it preserves the invariance of $r_{\mathcal{G}}(p)$, i.e.,
	\begin{eqnarray}
		\label{eq:1}
		\dot{r}_{\mathcal{G}}(p) = \frac{\partial r_{\mathcal{G}}(p)}{\partial p}v = 0.
	\end{eqnarray}
	A motion is said to be \textit{trivial} if it preserves rigidity for any framework. For example, if each entry of $r_\mathcal{G}(p)$ is a distance, i.e., $\|p_i-p_j\|$, then the trivial motions consist of uniform rotations and translations of the whole framework  \cite{hendrickson1992conditions}.
	
	\begin{definition} \label{definition IR}
		In $\mathbb{R}^d$, a framework $(\mathcal{G}, p)$ is said to be \textit{infinitesimally rigid} if all of its infinitesimal motions are trivial.
	\end{definition}	
	
	The \textit{rigidity matrix} is defined as $R(p)\triangleq\frac{\partial r_{\mathcal{G}}(p)}{\partial p}$. Then, (\ref{eq:1}) can be rewritten as $\dot{r}_{\mathcal{G}}(p) = R(p)\dot{p} = 0$. It can be obtained that a framework $(\mathcal{G}, p)$ in $\mathbb{R}^d$ is infinitesimally rigid if and only if $\rank(R(p)) = dn-T$, where $T$ is the dimension of the space spanned by trivial motions.
	
	So far, we have introduced some relevant rigidity concepts based on a given rigidity function. Throughout this paper, we mainly focus on the case in $\mathbb{R}^2$. In following subsections, we will introduce two specific rigidity functions, i.e., bearing rigidity function and signed angle rigidity function, based on which similar rigidity concepts can be induced naturally.

 	\begin{figure}
		\centering
		\includegraphics[width=.6\linewidth]{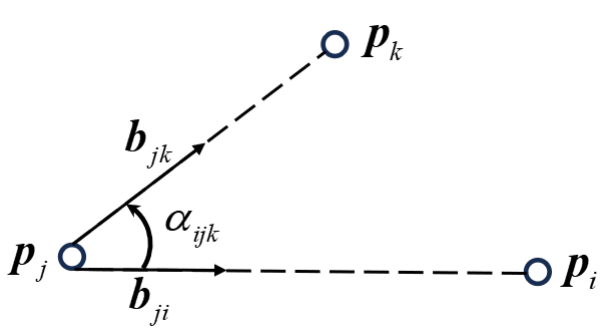} 
		\caption{{\footnotesize Example of the bearing and the signed angle among three agents. }}\label{fig: example of signed angle and bearing}
	\end{figure}
	
	\subsection{Bearing Rigidity}
	
	Given a framework $(\mathcal{G}, p)$ in $\mathbb{R}^2$, according to \cite{zhao2015bearing}, the \textit{bearing rigidity function} is described as:
	\begin{eqnarray}
		\label{eq:2}
		B_{\mathcal{G}}(p) = \left(\cdots, b^{\top}_{ji}, \cdots\right)^{\top}, (j,i)\in\mathcal{E},
	\end{eqnarray}
	where $b_{ji} \triangleq b_{ji}(p)  = \frac{p_i-p_j}{\|p_i-p_j\|} \in \mathbb{R}^{2}$ is the bearing vector from $p_j$ to $p_i$ (see Fig. \ref{fig: example of signed angle and bearing}). The {\it bearing rigidity matrix} is
	\begin{eqnarray}
		R_B(p) \triangleq  {\rm diag}\left(\frac{P_{ji}}{\|e_{ji}\|}\right) \bar{H} \in \mathbb{R}^{2|\mathcal{E}|\times 2n},
	\end{eqnarray}
	where $e_{ji} \triangleq p_i-p_j$, $P_{ji} \triangleq P(b_{ji}) = I_2-b_{ji}b^{\top}_{ji}$ is an orthogonal projection matrix projecting any vector onto the orthogonal space of $b_{ji}$, $ \bar{H} \triangleq \left(H \otimes I_2\right)$, $H = \left[h_{ji}\right]$ is the incidence matrix whose rows are indexed by edges of $\mathcal{G}$ and columns are indexed by vertices of $\mathcal{G}$ with an orientation:
	$h_{ji} = 1$ if vertex $i$ is the end of the $j$th edge, $h_{ji} = -1$ if vertex $i$ is the start of the $j$th edge, and $h_{ji} = 0$ otherwise. 
	
	Since bearing-preserving motions of a framework in $\mathbb{R}^2$ include 2-dimensional uniform translations and 1-dimensional uniform scalings \cite{zhao2015bearing}, a framework is infinitesimally bearing rigid if and only if $\rank(R_B(p)) = 2n-3$.

	\subsection{Signed Angle Rigidity} 
	Given a framework $(\mathcal{G}, p)$ in $\mathbb{R}^2$, according to \cite{chen2020angle}, the \textit{signed angle rigidity function} is defined as:
	\begin{eqnarray}
		\label{eq:4}
		S_{\mathcal{G}}(p) = \left(\cdots, \alpha_{ijk}, \cdots\right)^\top, (i,j,k) \in \mathcal{T}_{\mathcal{G}},
	\end{eqnarray}
	where $\alpha_{ijk} \triangleq \alpha_{ijk}(p) \in [0, 2\pi) $ describes the signed angle from $b_{ji}$ to $b_{jk}$ in the counter-clockwise direction (see Fig. \ref{fig: example of signed angle and bearing}), and is computed by
	\begin{eqnarray}
		\label{eq:signed angle}
		\alpha_{ijk} & =\left\{\begin{array}{ll}
			{\rm arccos}\left(b^{\top}_{jk}b_{ji}\right) & {\rm if}~ b^{\top}_{jk}b^{\bot}_{ji} \geq 0, \\
			2\pi-{\rm arccos}\left(b^{\top}_{jk}b_{ji}\right) & {\rm otherwise},
		\end{array}\right.
	\end{eqnarray} 
	$\mathcal{T}_{\mathcal{G}} = \left\{(i,j,k)\in\mathcal{V}^3: (j,i),(j,k)\in\mathcal{E}, i<k\right\}$ is a set indexing all signed angles in $(\mathcal{G}, p)$, here $(i,j,k)$ corresponds to the signed angle $\alpha_{ijk}$. By differentiating $S_\mathcal{G}(p)$ w.r.t. $p$, one has the {\it signed angle rigidity matrix}:
	\begin{eqnarray} \label{signed angle rigid matrix}
		R^{\mathcal{T}_{\mathcal{G}}}_S(p) \triangleq {\rm diag}\left(-\frac{1}{{\rm sin}\alpha_{ijk}}\right)R_b(p)R_B(p)\in \mathbb{R}^{|\mathcal{T}_{\mathcal{G}}|\times 2n},
	\end{eqnarray}
	where $R_B(p)$ is the bearing rigidity matrix and $R_b(p) \in \mathbb{R}^{|\mathcal{T}_{\mathcal{G}}| \times 2|\mathcal{E}|}$ can be written as
	\begin{eqnarray}
		\bordermatrix{%
			&\cdots & b_{ji}  &\cdots& b_{jk}    &\cdots \cr
			\hspace{5mm}    \cdots &\cdots   &\cdots &\cdots  &\cdots  & \cdots\cr 
			\hspace{5mm}	\cdots  &\cdots   &\cdots &\cdots  &\cdots  & \cdots\cr
			\hspace{3mm} \alpha_{ijk} & \boldsymbol{0}    &b_{jk}^\top   &\boldsymbol{0}  &b_{ji}^\top & \boldsymbol{0} \cr
			\hspace{5mm}	\cdots    &\cdots   &\cdots &\cdots  &\cdots  & \cdots\cr
			\hspace{5mm}    \cdots & \cdots &\cdots &\cdots  &\cdots  & \cdots
		}\notag
	\end{eqnarray}
	whose rows are indexed by $\mathcal{T}_{\mathcal{G}}$ and columns are indexed by the edges of $\mathcal{G}$. 
	
	It should be noted that the existence of $\sin \alpha_{ijk} = 0$ in (\ref{signed angle rigid matrix}) seems to render the signed angle rigidity matrix $R^{\mathcal{T}_{\mathcal{G}}}_S(p)$ invalid. In order to eliminate this singularity, we derive another expression for the signed angle rigidity matrix.

    \begin{lemma} \label{lemma: new form of signed angle rigidity matrix}
        The signed angle rigidity matrix (\ref{signed angle rigid matrix}) can be expressed as 
        \begin{eqnarray} \label{signed angle rigidity matrix 2}
            R^{\mathcal{T}_{\mathcal{G}}}_S(p) = \bar{R}_{\mathcal{T}_{\mathcal{G}}}(p)\bar{H},
        \end{eqnarray}
        where $\bar{R}_{\mathcal{T}_{\mathcal{G}}}(p) \in \mathbb{R}^{|\mathcal{T}_{\mathcal{G}}|\times 2|\mathcal{E}|}$ is
	\begin{eqnarray}
		\bordermatrix{%
			&\cdots & b_{ji}  &\cdots& b_{jk}    &\cdots \cr
			\hspace{5mm}    \cdots &\cdots   &\cdots &\cdots  &\cdots  & \cdots\cr 
			\hspace{5mm}	\cdots  &\cdots   &\cdots &\cdots  &\cdots  & \cdots\cr
			\hspace{3mm} \alpha_{ijk} & \boldsymbol{0}    &\frac{b_{ji}^\top\mathscr{R}_o(\frac{\pi}{2})  }{\|e_{ji}\|} &\boldsymbol{0}  &-\frac{b_{jk}^\top\mathscr{R}_o(\frac{\pi}{2})}{\|e_{jk}\|} & \boldsymbol{0} \cr
			\hspace{5mm}	\cdots    &\cdots   &\cdots &\cdots  &\cdots  & \cdots\cr
			\hspace{5mm}    \cdots & \cdots &\cdots &\cdots  &\cdots  & \cdots
		}.\notag
	\end{eqnarray}
    \end{lemma}
    \begin{IEEEproof}
        See Appendix \ref{Proof of Lemma lemma: new form of signed angle rigidity matrix}.
    \end{IEEEproof}
    
    In bearing rigidity theory, a framework is infinitesimally rigid if and only if the rank of the related rigidity matrix is $2n-3$. One may question if there is a similar conclusion in signed angle rigidity theory. In \cite{chen2020angle}, the authors proposed the signed angle rigidity theory for angularities and showed that an angularity is infinitesimally signed angle rigid if and only if $\rank(R^{\bar{\mathcal{T}}_{\mathcal{G}}}_S(p)) = 2n-4$, where an angularity involves a pre-designed subset $\bar{\mathcal{T}}_{\mathcal{G}} \subseteq \mathcal{T}_{\mathcal{G}}$ specifying constrained signed angles with the assumption that $\sin\alpha_{ijk} \neq 0$ for $(i,j,k) \in \bar{\mathcal{T}}_{\mathcal{G}}$. 
    
    In order to build connections with bearing rigidity theory in a unified architecture, we utilize all the signed angles in the framework to study the signed angle rigidity property without imposing any assumptions on these signed angle constraints. That is, given a framework $(\mathcal{G}, p)$, it becomes certain whether there exists a subset $\bar{\mathcal{T}}_{\mathcal{G}} \subseteq \mathcal{T}_{\mathcal{G}}$ that validates (infinitesimal, global) signed angle rigidity.
    
     To better discuss the existence of such a subset, we will propose several definitions in the following subsection. 
	
	\subsection{Angle Index Set and Angle Index Graph}
	
	Given a graph $\mathcal{G}$, we refer to a subset $\bar{\mathcal{T}}_{\mathcal{G}} \subseteq \mathcal{T}_{\mathcal{G}}$ as the \textit{angle index set} (AIS) and present the following definitions.
 
	\begin{definition} \label{definition: GAISs}
		For a given framework $(\mathcal{G},p)$, an AIS $\bar{\mathcal{T}}_{\mathcal{G}} \subseteq \mathcal{T}_{\mathcal{G}}$ is said to be the \textit{global angle index set} (GAIS) if for any $q\in \mathbb{R}^{2n}$ such that $\alpha_{ijk}(q) = \alpha_{ijk}(p), (i,j,k) \in \bar{\mathcal{T}}_{\mathcal{G}}$, it always holds that $S_{\mathcal{G}}(q) = S_{\mathcal{G}}(p)$.
	\end{definition}
        
        By definition, $\mathcal{T}_{\mathcal{G}}$ is always a GAIS. 
        
	\begin{definition} \label{definition: RAISs}
		For a given framework $(\mathcal{G},p)$, an AIS $\mathcal{T}^*_{\mathcal{G}} \subseteq \mathcal{T}_{\mathcal{G}}$ is said to be a \textit{regular angle index set} (RAIS) if $\rank(R^{\mathcal{T}^*_{\mathcal{G}}}_S(p)) =  \rank(R_S^{\mathcal{T}_\mathcal{G}}(p))$.
	\end{definition}

    \begin{definition}
    Given an RAIS $\mathcal{T}^*_{\mathcal{G}}$, if there does not exist another RAIS with fewer elements than $\mathcal{T}^*_{\mathcal{G}}$, then it is said to be a {\it minimal} RAIS. The same goes for GAIS.
    \end{definition}

    Since an angular constraint is related to two edges in the framework, it is necessary to evaluate an AIS from the perspective of the edge set. We use the Hash function $H(i,j,|\mathcal{V}|) = (\min(\{i,j\})-1)\times |\mathcal{V}|+\max(\{i,j\})$ to encode the edge $(i,j)$ of a given graph $\mathcal{G}$ sequentially, and propose the following definitions
	
    \begin{definition}[Angle Index Graph]\label{de angle index graph}
		Given a graph $\mathcal{G}$ and an AIS $\bar{\mathcal{T}}_{\mathcal{G}}\subseteq \mathcal{T}_{\mathcal{G}}$, the undirected {\it angle index graph} is denoted by $\mathcal{G}_A(\bar{\mathcal{T}}_{\mathcal{G}}) = \left(\mathcal{V}_A(\mathcal{G}), \mathcal{E}_A(\bar{\mathcal{T}}_{\mathcal{G}}) \right)$, where the vertex set is $\mathcal{V}_A(\mathcal{G}) = \{a_{ij} = H(i,j, |\mathcal{V}|): (i,j)\in\mathcal{E}\}$, $a_{ij}\in\mathcal{V}_A(\mathcal{G})$ is equivalent to $a_{ji}\in\mathcal{V}_A(\mathcal{G})$, and the edge set is $\mathcal{E}_A(\bar{\mathcal{T}}_{\mathcal{G}}) = \{(a_{ji},a_{jk}):a_{ji},a_{jk}\in\mathcal{V}_A(\mathcal{G}),(i,j,k) \in \bar{\mathcal{T}}_{\mathcal{G}}\}$.
    \end{definition}

    If $\mathcal{G}$ and $\bar{\mathcal{T}}_{\mathcal{G}}$ are clear from the context, we simply denote $\mathcal{G}_A(\bar{\mathcal{T}}_{\mathcal{G}}) = \left(\mathcal{V}_A(\mathcal{G}), \mathcal{E}_A(\bar{\mathcal{T}}_{\mathcal{G}}) \right)$ as $\mathcal{G}_A=\left(\mathcal{V}_A, \mathcal{E}_A \right)$.
    
    \begin{definition}\label{de conected AIS} Given a graph $\mathcal{G}$ and an AIS $\bar{\mathcal{T}}_{\mathcal{G}}\subseteq \mathcal{T}_{\mathcal{G}}$, $\bar{\mathcal{T}}_{\mathcal{G}}$ is said to be \textit{angle connected} if the angle index graph $\left(\mathcal{V}_A(\mathcal{G}), \mathcal{E}_A(\bar{\mathcal{T}}_{\mathcal{G}}) \right)$ is connected.
	\end{definition}
    
	\begin{figure}
		\centering
		\includegraphics[width=.9\linewidth]{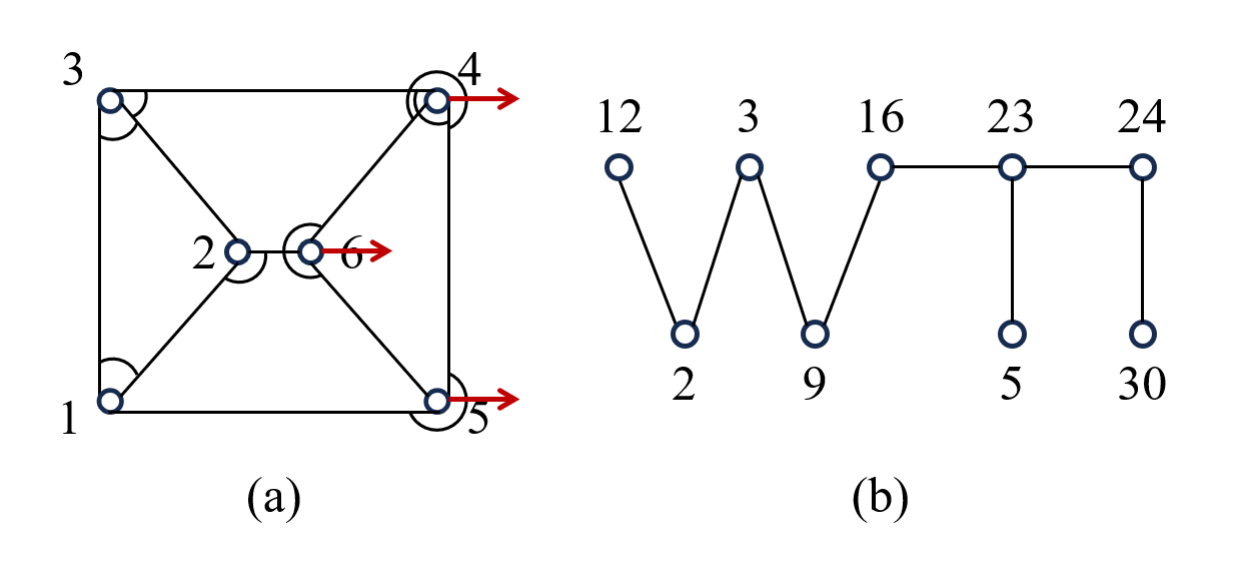} 
		\caption{{\footnotesize (a) A non-infinitesimally signed angle rigid framework $(\mathcal{G},p)$ and an angle connected set $\bar{\mathcal{T}}_{\mathcal{G}} = \{(2,1,3), (1,2,6), (1,3,2), (2,3,4), (3,4,5), (\\5,4,6), (1,5,4), (4,6,5) \}$, where the red arrows stand for non-trivial infinitesimal signed angle motions. (b) The angle index graph of $\mathcal{G}$ and $\bar{\mathcal{T}}_{\mathcal{G}}$, where vertex $12$ satisfies $ 12 = a_{26} = H(2,6,6) = (2-1)\times6+6$. }}
        \label{fig: example of angle index graph}
	\end{figure}
 
    Note that $\mathcal{G}_A(\bar{\mathcal{T}}_{\mathcal{G}}) = \left(\mathcal{V}_A(\mathcal{G}), \mathcal{E}_A(\bar{\mathcal{T}}_{\mathcal{G}}) \right)$ is formed by mapping the edge $(i,j) \in \mathcal{E}$ to the vertex in $\mathcal{G}_A(\bar{\mathcal{T}}_{\mathcal{G}})$ and the index $(i,j,k) \in \bar{\mathcal{T}}_{\mathcal{G}}$ to the undirected edge in $\mathcal{G}_A(\bar{\mathcal{T}}_{\mathcal{G}})$, respectively. Furthermore, for a given graph $\mathcal{G}$, an angle index graph $\mathcal{G}_A(\bar{\mathcal{T}}_{\mathcal{G}})$ can uniquely deduce an AIS. Fig. \ref{fig: example of angle index graph} illustrates an example of the angle connected set $\bar{\mathcal{T}}_{\mathcal{G}} \subseteq \mathcal{T}_{\mathcal{G}}$, where the angular constraints between edges and their interrelationships are visually presented.
 
	 Throughout this paper, we use shorthands SAR, GSAR, and ISAR for signed angle rigid (signed angle rigidity), globally signed angle rigid (global signed angle rigidity), and infinitesimally signed angle rigid (infinitesimal signed angle rigidity), respectively. Similarly, BR, GBR, and IBR are shorthands for corresponding concepts in bearing rigidity theory.
	
	\section{Graphs for Signed Angle Rigidity}\label{sec:3}
	
	In this section, we establish the equivalence between signed angle constraints and bearing constraints in the sense that a framework $(\mathcal{G}, p)$ is IBR, GBR, and BR if and only if it is ISAR, GSAR, and SAR, respectively. Based on these relationships, we derive necessary and sufficient graphical conditions for signed angle rigidity.
	
	\subsection{Equivalence on IBR and ISAR}
        In \cite{chen2020angle}, the authors showed that if $\sin \alpha_{ijk} \neq 0, (i,j,k) \in\mathcal{T}_{\mathcal{G}}$, the null space of the signed angle rigidity matrix $R^{\mathcal{T}_{\mathcal{G}}}_{S}(p)$ contains the following three linear spaces:
        \begin{eqnarray}
            \mathcal{M}_s &\hspace{-2mm}\triangleq\hspace{-2mm}& {\rm span} \{p\}, \\
            \mathcal{M}_r &\hspace{-2mm}\triangleq\hspace{-2mm}& {\rm span} \left\{(I_n \otimes \mathscr{R}_o(\frac{\pi}{2}))p\right\},\\
            \mathcal{M}_t &\hspace{-2mm}\triangleq\hspace{-2mm}& {\rm span}\left\{\textbf{1}_n\otimes(1,0)^\top,\textbf{1}_n \otimes (0,1)^\top\right\}.
        \end{eqnarray}
        According to Lemma \ref{lemma: new form of signed angle rigidity matrix}, it can be verified that the result remains valid even when there exist triplets $(i,j,k)\in \mathcal{T}_{\mathcal{G}}$ satisfying $\sin \alpha_{ijk} = 0$. Therefore, the following lemma can be derived from Definition \ref{definition IR} naturally, which implies that the trivial infinitesimal signed angle motions of a framework in $\mathbb{R}^2$ include 2-dimensional uniform translations, 1-dimensional uniform scalings, and 1-dimensional uniform rotations.

        \begin{lemma} \label{lemma: trivial motions of signed angle rigidity 2}
            In $\mathbb{R}^2$, a framework $(\mathcal{G},p)$ is ISAR if and only if $\Null(R^{\mathcal{T}_{\mathcal{G}}}_S(p)) = \mathcal{M}_s\bigcup \mathcal{M}_r \bigcup \mathcal{M}_t$.
        \end{lemma}

    As shown in Fig. \ref{fig: example of angle index graph}(a), the framework is not ISAR since there exists a non-trivial infinitesimal signed angle motion. Lemma \ref{lemma: trivial motions of signed angle rigidity 2} also indicates that a framework $(\mathcal{G},p)$ in $\mathbb{R}^2$ is ISAR if and only if $\rank(R^{\mathcal{T}_{\mathcal{G}}}_S(p))$ = $2n-4$. Next, we present a key property of ISAR frameworks.
        \begin{theorem}  \label{corollary ibr isar and iar}
		In $\mathbb{R}^2$, a framework $(\mathcal{G},p)$ is ISAR if and only if it is IBR.
	\end{theorem}
        \begin{IEEEproof}
            See Appendix \ref{Proof of Theorem 1}.
        \end{IEEEproof}
    
	Theorem \ref{corollary ibr isar and iar} implies that the conditions on the framework for bearing constraints and signed angle constraints to exclude non-trivial infinitesimal motions are the same. In \cite{zhao2015bearing}, it has been shown that the shape of an IBR framework can always be uniquely fixed by all bearings among neighboring nodes. However, Theorem \ref{corollary ibr isar and iar} does not imply that an ISAR framework can fix its shape by using all its angular information. To seek the conditions on frameworks for signed angles to constrain a shape uniquely, we study the relationship between GSAR and GBR in the next subsection.
	
	\subsection{Equivalence on GSAR and GBR}
	In this subsection, we establish the equivalence among GSAR, GBR, SAR and BR. Before that, we propose a graphical condition, under which signed angle constraints are equivalent to bearing constraints regardless of uniform rotations.
 
	\begin{lemma} \label{lemma equivalence between sar and br}
		Given a connected graph $\mathcal{G}$ and two configurations $p, q \in \mathbb{R}^{2n}$, $S_{\mathcal{G}}(q) = S_{\mathcal{G}}(p)$ if and only if $B_{\mathcal{G}}(q) = \left(I_n \otimes \mathscr{R}_o(\theta) \right) B_{\mathcal{G}}(p)$ for some $\theta \in [0,2\pi)$.
	\end{lemma}
	
	\begin{IEEEproof}		
		See Appendix \ref{proof of Lemma 3}.
	\end{IEEEproof}
	
	Lemma \ref{lemma equivalence between sar and br} implies that, despite being measured in different local coordinate frames, signed angle constraints are able to function as bearing constraints in a certain coordinate frame as long as the graph is connected. This is an interesting property with great potential in coordinate-free formation control and network localization, but surprisingly, not found in the literature. 
	
	Now we are ready to show the equivalence between GBR and GSAR.
	
	\begin{theorem} \label{theorem equivalence between BR and SAR}
		A framework $(\mathcal{G}, p)$ in $\mathbb{R}^2$ is BR (GBR) if and only if it is SAR (GSAR).
	\end{theorem}
	
	\begin{IEEEproof}
		Sufficiency. Suppose $(\mathcal{G}, p)$ is SAR. Then, there exists a neighborhood $U_p$ of $p$ such that for any $q \in S^{-1}_{\mathcal{G}}\left(S_{\mathcal{G}}(p)\right)\bigcap U_p$, it always holds $S_{\mathcal{K}}(q) = S_{\mathcal{K}}(p)$. It suffices to show that for this $U_p$, consider any $q \in B^{-1}_{\mathcal{G}}(B_{\mathcal{G}}(p))\bigcap U_p$, it always holds $B_{\mathcal{K}}(q) = B_{\mathcal{K}}(p)$. 
		
		From $B_{\mathcal{G}}\left( q \right) = B_{\mathcal{G}}(p)$ and (\ref{eq:signed angle}), we have $S_{\mathcal{G}}(q) = S_{\mathcal{G}}(p)$, which also implies $S_{\mathcal{K}}(q) = S_{\mathcal{K}}(p)$. Since $b_{ji}(p) = b_{ji}(q), (j,i) \in \mathcal{E}$ and all signed angles of $\mathcal{K}(p)$ are determined, we can conclude that $B_{\mathcal{K}}\left( q \right) = B_{\mathcal{K}}(p)$.
		
		Necessity. Suppose $(\mathcal{G}, p)$ is BR. Then, there exists a neighborhood $U_p$ of $p$ such that for any $q \in B^{-1}_{\mathcal{G}}\left(B_{\mathcal{G}}(p)\right)\bigcap U_p$, it always holds $B_{\mathcal{K}}(q) = B_{\mathcal{K}}(p)$. It suffices to show that for this $U_p$, consider any $q \in S^{-1}_{\mathcal{G}}(S_{\mathcal{G}}(p))\bigcap U_p$, it always holds $S_{\mathcal{K}}(q) = S_{\mathcal{K}}(p)$. Note that the embedded graph of the BR framework must be connected, then it follows from Lemma \ref{lemma equivalence between sar and br} and $S_{\mathcal{G}}(q) = S_{\mathcal{G}}(p)$ that $B_{\mathcal{G}}\left( q \right) = \left(I_n \otimes \mathscr{R}_o(\theta) \right)B_{\mathcal{G}}(p)$ for some $\theta \in [0,2\pi)$. Since the framework is BR, it follows that $B_{\mathcal{K}}\left( q \right) = \left(I_n \otimes \mathscr{R}_o(\theta) \right)B_{\mathcal{K}}(p)$. Thus, for any $i,j,k\in\mathcal{V}$,
        \begin{eqnarray} \label{eq:bearin to signed angle1}
			b^{\top}_{jk}(q)b_{ji}(q) \hspace{-2mm}&=&\hspace{-2mm} b^{\top}_{jk}(p)\mathscr{R}^{\top}_o(\theta) \mathscr{R}_o(\theta)b_{ji}(p) \notag\\
            \hspace{-2mm}&=&\hspace{-2mm} b^{\top}_{jk}(p)b_{ji}(p),
		\end{eqnarray}
        and
		\begin{eqnarray} \label{eq:bearin to signed angle2}
			b^{\top}_{jk}(q)\mathscr{R}_o(\frac{\pi}{2})b_{ji}(q) \hspace{-2mm}&=&\hspace{-2mm} b^{\top}_{jk}(p)\mathscr{R}^{\top}_o(\theta)\mathscr{R}_o(\frac{\pi}{2})\mathscr{R}_o(\theta)b_{ji}(p) \notag\\
			&=&\hspace{-2mm} b^{\top}_{jk}(p)\mathscr{R}_o(\frac{\pi}{2})b_{ji}(p).
		\end{eqnarray}
		According to (\ref{eq:signed angle}), we further have $S_{\mathcal{K}}(q) = S_{\mathcal{K}}(p)$.
		
		Similar to the proof of the equivalence between BR and SAR, we can derive the same result reflected on GBR and GSAR by extending $U_p$ to $\mathbb{R}^{2n}$. 
	\end{IEEEproof}
 
	Combining Theorems \ref{corollary ibr isar and iar} and \ref{theorem equivalence between BR and SAR}, together with \cite[Theorem 3]{zhao2015bearing} and \cite[Theorem 5]{zhao2015bearing}, the following results hold.

	\begin{corollary} \label{corollary relationship between isar and gsar}
		If a framework $(\mathcal{G}, p)$ in $\mathbb{R}^2$ is ISAR, then it is GSAR.
	\end{corollary}
    
	\begin{corollary} \label{corollary equivalence between gsar and sar}
		A framework $(\mathcal{G}, p)$ in $\mathbb{R}^2$ is GSAR if and only if it is SAR.
	\end{corollary}
 
\begin{figure}[hthb]
		\centering
		\includegraphics[width=.9\linewidth]{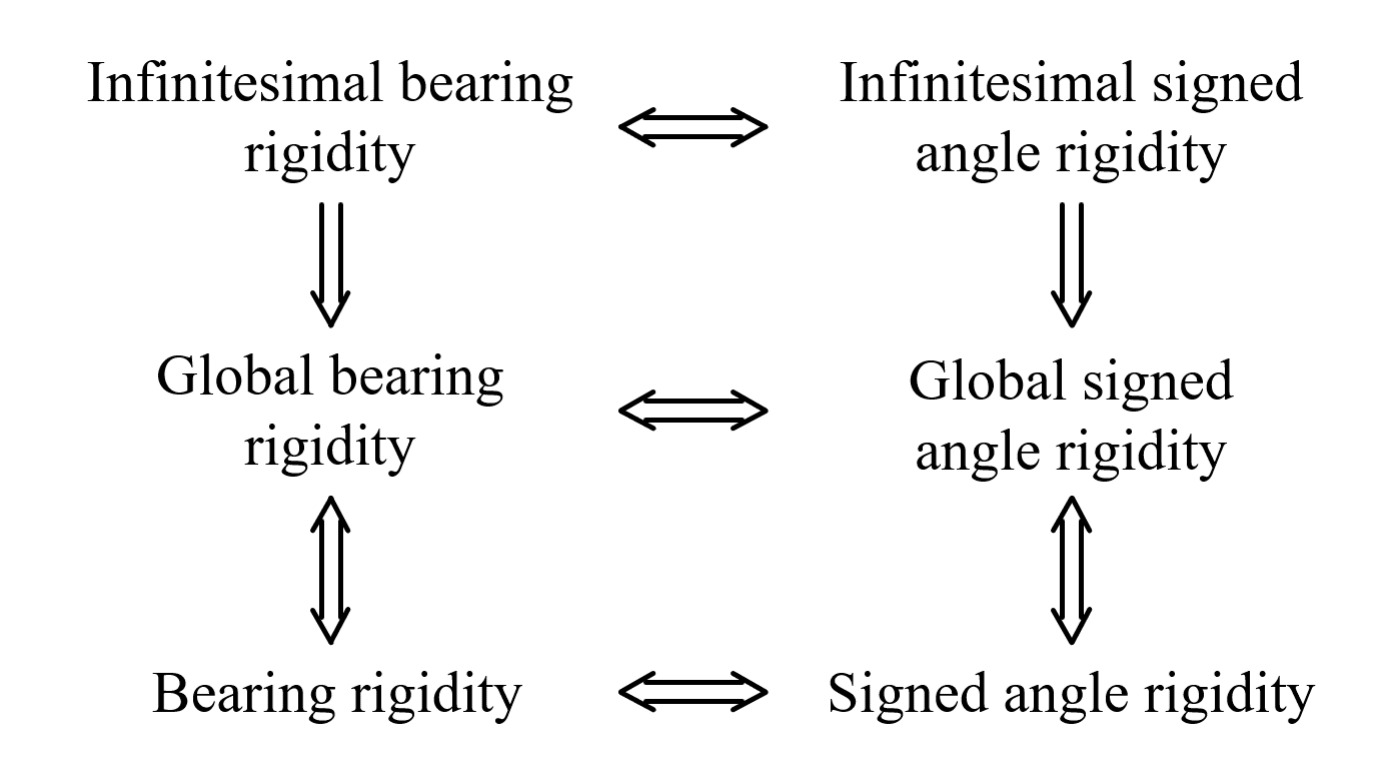}
		\caption{A relation conclusion graph of different conceptions of bearing rigidity and signed angle rigidity in 2-D}
		\label{relaiton conclusion graph}
		\centering
	\end{figure}

        To clearly illustrate the results stated above, we summarize the relationships between different concepts of bearing rigidity and signed angle rigidity in Fig. \ref{relaiton conclusion graph}.
 	
        \begin{remark}
          Under the context of angularity, the authors in \cite{chen2020angle} and \cite{zhou2006combinatorial} proved that signed angle rigidity is not equivalent to global signed angle rigidity, and infinitesimal signed rigidity does not imply global signed angle rigidity. In this paper, we arrive at different conclusions by taking into account all the signed angles in the frameworks.
        \end{remark}
	
	\begin{remark}
		We have established complete equivalence on different rigidity properties based on signed angle constraints and bearing constraints. One may question if such equivalence holds for unsigned angles in \cite{jing2019angle}. Similarly to the conclusion of shape fixability in bearing rigidity theory \cite{zhao2015bearing}, it has been proven in \cite{jing2021angle} that a framework in $\mathbb{R}^2$ can be uniquely determined by unsigned angles if and only if it is infinitesimally angle rigid and globally angle rigid. Although several illustrations induce the conjecture that infinitesimal angle rigidity may imply global angle rigidity (like ISAR implie GSAR), this is yet to be proven.
	\end{remark}
	
	\subsection{The Uniqueness of Shape}
	
	Based on the above results, we next answer an important question: under what conditions can the shape of a framework be uniquely determined by signed angle constraints?
	
	Given a configuration $p \in \mathbb{R}^{2n}$, the set of configurations that form the same shape as $p$ can be described as:
	\begin{multline}
		\mathscr{E}(p) = \{q\in\mathbb{R}^{2n}:q = c(I_n \otimes\mathscr{R}_o(\theta))p+\textbf{1}_n \otimes \xi, \theta \in [0,2\pi), \\
		c\in \mathbb{R}\backslash \{0\}, \xi \in \mathbb{R}^{2} \}.
	\end{multline}
	If $p\in\mathbb{R}^{2n}$ is non-degenerate, then $\mathscr{E}(p)$ is a 4-dimensional smooth manifold. Given a graph $\mathcal{G}$ with $n$ vertices, the configuration set of frameworks with the same signed angles as $(\mathcal{G}, p)$ can be depicted as:
        \begin{eqnarray}
            \mathscr{E}_{\mathcal{G}}(p) = \{q\in\mathbb{R}^{2n}:\alpha_{ijk}(q) = \alpha_{ijk}(p), (i,j,k)\in \mathcal{T}_{\mathcal{G}}\}.
        \end{eqnarray}
	Then, we have the following lemma naturally.
	
	\begin{lemma} \label{lemma shape fix manifold}
		The shape of a framework $(\mathcal{G},p)$ in $\mathbb{R}^{2}$ can be uniquely determined by signed angles up to uniform rotations, translations, and scalings if and only if $\mathscr{E}(p) = \mathscr{E}_{\mathcal{G}}(p)$.
	\end{lemma}

        The following theorem presents the necessary and sufficient conditions for the shape fixability of frameworks under signed angle constraints, from the perspectives of ISAR and GSAR, respectively.
	
	\begin{theorem} \label{theorem shape fix}
		Given a framework $(\mathcal{G},p)$ in $\mathbb{R}^{2}$, the following statements are equivalent:
		
		\begin{itemize}
			\item [(\rmnum{1})] $(\mathcal{G},p)$ is ISAR; 
			\item [(\rmnum{2})] $(\mathcal{G},p)$ is non-degenerate and GSAR;
			\item [(\rmnum{3})] the shape of $(\mathcal{G},p)$ can be uniquely determined by signed angles up to uniform rotations, translations, and scalings.
		\end{itemize}
	\end{theorem}
	
	\begin{IEEEproof}
		See Appendix \ref{proof of Theorem 3}.
	\end{IEEEproof}
	
	Given a framework $(\mathcal{G},p)$ in $\mathbb{R}^2$, Theorem \ref{theorem shape fix} suggests that: (\rmnum{1}) ISAR, which can be examined by computing the rank of the signed angle rigidity matrix, is able to determine whether the shape of a framework can be specified by signed angles uniquely; (\rmnum{2}) ISAR (IBR) and GSAR (GBR) imply each other when $p$ is non-degenerate, which is interesting but not found in the literature.

    However, given the above results, computing the rank of the rigidity matrix for large-scale networks may be computationally prohibitive. In the next subsection, we study how to efficiently construct graphs that are ISAR (GSAR, SAR) for almost all configurations.

	\subsection{Signed Angle Rigid Graphs}
        To begin with, we formally introduce the definition of signed angle rigidity for graphs.
        \begin{definition}
            A graph is said to be SAR (ISAR, GSAR) if there exists at least one configuration $p\in\mathbb{R}^{2n}$ such that $(\mathcal{G},p)$ is SAR (ISAR, GSAR).
        \end{definition}

        In \cite[Lemma 2]{zhao2017laman}, IBR has been shown as a generic property of graphs. More specifically, if $(\mathcal{G},p)$ is IBR for a configuration $p\in\mathbb{R}^{2n}$, then $(\mathcal{G},q)$ is IBR for almost all $q\in\mathbb{R}^{2n}$. According to \cite{zhao2015bearing}, such a graphic property holds for GBR and BR as well. From the relationship between bearing rigidity and signed angle rigidity shown in Fig. \ref{relaiton conclusion graph}, it can be observed that ISAR, GSAR, and SAR are essentially determined by the graph rather than the configuration, thus are generic properties of graphs. We summarize this result in the following theorem.
        
	\begin{theorem} \label{theorem: generic graph property of a graph}
		If a graph $\mathcal{G}$ is ISAR (GSAR, SAR), then the framework $(\mathcal{G},p)$ is ISAR (GSAR, SAR) for almost all configurations $p \in \mathbb{R}^{2n}$.
	\end{theorem}

        The equivalent graphical description for bearing rigidity has been well understood in \cite{zhao2017laman, trinh2019minimal}. Inspired by this, we derive the equivalent graphical condition for signed angle rigidity. Before that, we introduce a class of well-known graphs.

	\begin{definition}[Laman Graph \cite{zhao2017laman,jackson2007notes}]
		A $n$-point $(n\geqslant3)$ graph $\mathcal{L}_n = (\mathcal{V}_n, \mathcal{E}_n)$ satisfying $|\mathcal{E}_n| = 2n - 3$ is a Laman graph if and only if $e \leqslant 2v - 3$ for every subgraph of $\mathcal{L}_n$ having $v$ vertices and $e$ edges.
	\end{definition}
	
	Combining Theorem \ref{corollary ibr isar and iar} and \cite[Theorem 2]{zhao2017laman}, we present the following result, which implies the shape of the graphs containing Laman spanning subgraphs can be determined by signed angle constraints for almost all configurations.
	\begin{theorem} \label{theorem: ISAR laman subframeworks}
		A graph $\mathcal{G}$ is ISAR if and only if it contains Laman spanning subgraphs.
	\end{theorem}

        In \cite{trinh2019minimal}, the authors have shown that minimally IBR graphs are Laman graphs. However, they did not give equivalent graphical conditions forminimally BR graphs. According to Theorem \ref{theorem: generic graph property of a graph} and Theorem \ref{theorem shape fix}, ISAR and GSAR imply each other for graphs. As a result, a graph $\mathcal{G}$ is GSAR if and only if it contains Laman spanning subgraphs. Based on the equivalence between bearing rigidity and signed angle rigidity, together with \cite[Theorem 4.1]{trinh2019minimal}, we have the following theorem.

        
        \begin{theorem} \label{theorem: MSAR is equicalent to Laman graph}
            A graph $\mathcal{G}$ is minimally SAR ( ISAR, GSAR, IBR, GBR, BR)  if and only if it is a Laman graph.
        \end{theorem}
        
        Now, it is clear which type of frameworks can be uniquely determined by signed angles. In order to exclude redundant angle constraints, the problem of finding a minimal set of signed angle constraints for shape determination will be investigated in the next section.

        \begin{remark} \label{remark :1}
            For any ISAR framework $(\mathcal{G},p)$ in $\mathbb{R}^2$, we can find an ISAR Laman spanning subframework in $(\mathcal{G},p)$ based on the distance rigidity matrix \cite{jackson2007notes}, where a framework is said to be the Laman framework if it is embedded by a Laman graph. Specifically, considering an ISAR framework $(\mathcal{G},p)$ in $\mathbb{R}^2$, according to \cite[Theorem 8]{zhao2015bearing} and Theorem \ref{corollary ibr isar and iar}, $(\mathcal{G},p)$ is infinitesimally distance rigid. Then there are $2n-3$ linearly independent rows in the distance rigidity matrix. According to \cite[Lemma 2.3]{jackson2007notes}, the edges corresponding to these $2n-3$ rows are also the edges of an ISAR Laman spanning subframework. 
        \end{remark}

	\section{Angle Index Sets for Signed Angle Rigidity} \label{sec:4}
        In this section, minimal GAISs and RAISs are investigated from the perspective of angle index graphs. Additionally, the relationship between AISs and signed angle rigidity is summarized as well.
	
    \subsection{Minimal GAISs for Signed Angle Rigidity}

        Given a framework $(\mathcal{G},p)$ in $\mathbb{R}^{2}$ and an AIS $\bar{\mathcal{T}}_{\mathcal{G}} \subseteq \mathcal{T}_{\mathcal{G}}$, the configuration set of frameworks with the same signed angles indexed by $\bar{\mathcal{T}}_{\mathcal{G}}$ as $(\mathcal{G}, p)$ can be depicted as:
        \begin{eqnarray}
            \mathscr{E}_{\bar{\mathcal{T}}_{\mathcal{G}}}(p) = \{q \in \mathbb{R}^{2n}:\alpha_{ijk}(q) = \alpha_{ijk}(p), (i,j,k) \in \bar{\mathcal{T}}_{\mathcal{G}}\},
        \end{eqnarray}
        From Definition \ref{definition: GAISs}, we have the following lemma naturally.

        \begin{lemma} \label{lemma: GAIS conditions}
            For a framework $(\mathcal{G},p)$ in $\mathbb{R}^2$, $\bar{\mathcal{T}}_{\mathcal{G}} \subseteq \mathcal{T}_{\mathcal{G}}$ is a GAIS if and only if it satisfies $\mathscr{E}_{\bar{\mathcal{T}}_{\mathcal{G}}}(p) = \mathscr{E}_{\mathcal{G}}(p)$.
        \end{lemma}

    \begin{figure}
		\centering
		\includegraphics[width=1\linewidth]{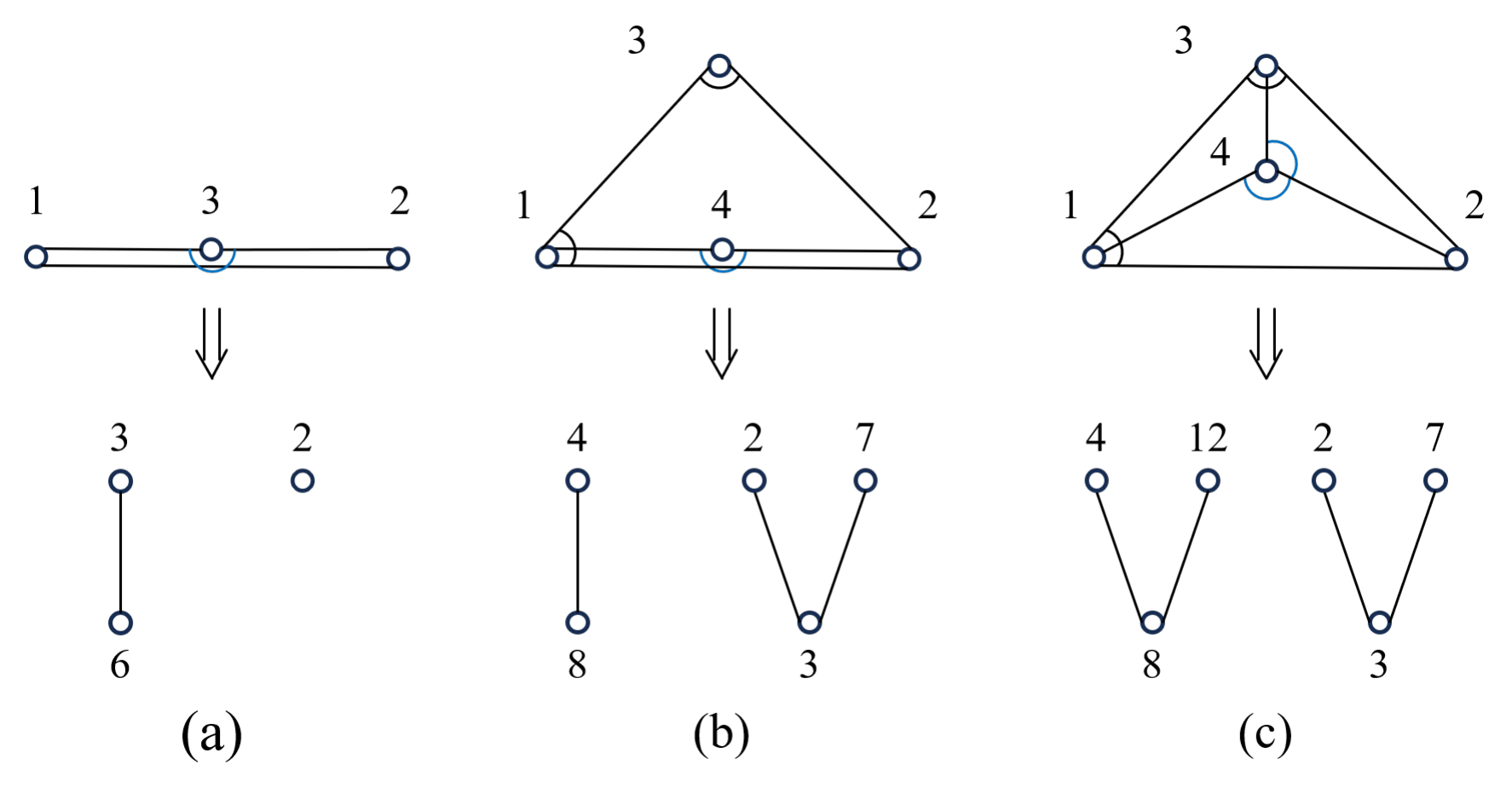}
		\caption{{\footnotesize Examples of GAISs whose corresponding angle index graphs are not angle connected. Black and blue are used to distinguish different components of $\bar{\mathcal{T}}_{\mathcal{G}}$. (a) A GSAR framework with a GAIS $\bar{\mathcal{T}}_{\mathcal{G}} = \left\{(1,3,2)\right\}$. (b) A non-GSAR framework with a GAIS $\bar{\mathcal{T}}_{\mathcal{G}} = \left\{(1,4,2), (1,3,2), (2,1,3)\right\}$. (c) A GSAR framework with a GAIS $\bar{\mathcal{T}}_{\mathcal{G}} = \{(1,4,2), (2,4,3), (1,3,2), (2,1,3)\}$.}}
		\label{fig: non-connected angle index graph}
	\end{figure} 

        The following result provides an intuitive condition to determine GAISs.
        
        \begin{lemma} \label{lemma: angle connected set is GAIS}
        For a framework $(\mathcal{G},p)$ in $\mathbb{R}^2$, an AIS $\bar{\mathcal{T}}_{\mathcal{G}}$ is a GAIS if it is angle connected.
        \end{lemma}

        \begin{IEEEproof}
            Since $\mathscr{E}_{\mathcal{G}}(p) \subseteq \mathscr{E}_{\bar{\mathcal{T}}_{\mathcal{G}}}(p)$, according to Lemma \ref{lemma: GAIS conditions}, it sufficient to show that for any $q\in\mathscr{E}_{\bar{\mathcal{T}}_{\mathcal{G}}}(p)$, one has $q \in \mathscr{E}_{\mathcal{G}}(p)$. 

            Note that for any edge $(a_{ji},a_{jk})$ of the angle index graph $(\mathcal{V}_A(\mathcal{G}), \mathcal{E}_A(\bar{\mathcal{T}}_{\mathcal{G}}))$, since $\alpha_{ijk}(q) = \alpha_{ijk}(p)$, one has $b_{ji}(q) = \mathscr{R}_o(\theta)b_{ji}(p)$ and $b_{jk}(q) = \mathscr{R}_o(\theta)b_{jk}(p)$ for some $\theta\in[0,2\pi)$. Since the AIS $\bar{\mathcal{T}}_{\mathcal{G}}$ is angle connected, it can be verified that $B_{\mathcal{G}}(q) = \left(I_n \otimes \mathscr{R}_o(\theta) \right) B_{\mathcal{G}}(p)$ for some $\theta\in[0,2\pi)$ by using similar techniques in the proof for necessity of Lemma \ref{lemma equivalence between sar and br}. As a result, we have $S_{\mathcal{G}}(q) = S_\mathcal{G}(p)$ and $q\in \mathscr{E}_{\mathcal{G}}(p)$.
        \end{IEEEproof}

        It is worth noting that the angle connectivity of GAISs is not necessary, see counterexamples in Fig. \ref{fig: non-connected angle index graph}. Interestingly, when the framework is embedded by a Laman graph (termed as a \textit{Laman framework}), one has a different conclusion.

        \begin{theorem} \label{lemma: AIS must be connected if and only if the graph is laman graph}
            For an ISAR Laman framework $(\mathcal{L}_n,p)$ in $\mathbb{R}^2$, $\bar{\mathcal{T}}_{\mathcal{L}_n}\subseteq \mathcal{T}_{\mathcal{L}_n}$ is a minimal GAIS if and only if $|\bar{\mathcal{T}}_{\mathcal{L}_n}| = 2n-4$ and $\bar{\mathcal{T}}_{\mathcal{L}_n}$ is angle connected.
        \end{theorem}
        \begin{IEEEproof}
            According to Theorem \ref{theorem shape fix} and Lemma \ref{lemma: GAIS conditions}, we know that for an ISAR framework $(\mathcal{G},p)$ in $\mathbb{R}^2$, $\bar{\mathcal{T}}_{\mathcal{G}} \subseteq \mathcal{T}_{\mathcal{G}}$ is a GAIS if and only if $\mathscr{E}_{\bar{\mathcal{T}}_{\mathcal{G}}}(p) = \mathscr{E}(p)$. Since $p$ is non-degenerate, $\mathscr{E}_{\bar{\mathcal{T}}_{\mathcal{G}}}(p)$ is a 4-dimensional smooth manifold. According to \cite[Proposition 3.10]{lee2013introduction}, the tangent space of $\mathscr{E}_{\bar{\mathcal{T}}_{\mathcal{G}}}(p)$ at $p$ is 4-dimensional, which indicates that $\rank (R^{\bar{\mathcal{T}}_{\mathcal{G}}}_S(p)) = 2n-4$. As a result, one has $|\bar{\mathcal{T}}_{\mathcal{G}}| \geq 2n-4$.
            
            The sufficiency follows from Lemma \ref{lemma: angle connected set is GAIS} and the fact that $|\bar{\mathcal{T}}_{\mathcal{L}_n}|$ is lower bounded by $2n-4$. To prove necessity, we denote the angle index graph of $\bar{\mathcal{T}}_{\mathcal{L}_n}$ and $\mathcal{L}_n = (\mathcal{V}_n,\mathcal{E}_n)$ by $\left(\mathcal{V}_A, \mathcal{E}_A \right)$ and show the following two cases contradict the fact that $\bar{\mathcal{T}}_{\mathcal{L}_n}$ is minimal to satisfy $\mathscr{E}_{\bar{\mathcal{T}}_{\mathcal{L}_n}}(p) = \mathscr{E}_{\mathcal{L}_n}(p)$:
            \begin{itemize}
                \item[(\rmnum{1})] $|\bar{\mathcal{T}}_{\mathcal{L}_n}| > 2n-4$ and $\bar{\mathcal{T}}_{\mathcal{L}_n}$ is angle connected;
                \item[(\rmnum{2})] $|\bar{\mathcal{T}}_{\mathcal{L}_n}| \geq 2n-4$ and $\bar{\mathcal{T}}_{\mathcal{L}_n}$ is not angle connected.
            \end{itemize}

            Case 1. Since $\bar{\mathcal{T}}_{\mathcal{L}_n}$ is angle connected, we can use the breadth first search method to find a spanning tree within $\left(\mathcal{V}_A, \mathcal{E}_A \right)$, it follows from Lemma \ref{lemma: angle connected set is GAIS} that the AIS $\bar{\mathcal{T}}^\dagger_{\mathcal{L}_n}$ corresponding to the spanning tree satisfies $\mathscr{E}_{\bar{\mathcal{T}}^\dagger_{\mathcal{L}_n}}(p) = \mathscr{E}_{\mathcal{L}_n}(p)$ and $|\bar{\mathcal{T}}^\dagger_{\mathcal{L}_n}| = |\mathcal{E}_n|-1$. Since the embedded graph is a Laman graph, one has $|\mathcal{E}_n| = 2n-3$, and consequently,
            \begin{eqnarray}
                |\bar{\mathcal{T}}^\dagger_{\mathcal{L}_n}| &\hspace{-2mm}=\hspace{-2mm}& |\mathcal{E}_n|-1 =  2n-4 < |\bar{\mathcal{T}}_{\mathcal{L}_n}|.
            \end{eqnarray}
            Hence, $\bar{\mathcal{T}}_{\mathcal{L}_n}$ is not minimal for $\mathscr{E}_{\bar{\mathcal{T}}_{\mathcal{L}_n}}(p) = \mathscr{E}_{\mathcal{L}_n}(p)$.
            
            Case 2. Without loss of generality, $\left(\mathcal{V}_A, \mathcal{E}_A \right)$ has only two connected components $\left(\mathcal{V}^1_A, \mathcal{E}^1_A \right)$ and $\left(\mathcal{V}^2_A, \mathcal{E}^2_A \right)$, where $\mathcal{V}^1_A\bigcup \mathcal{V}^2_A = \mathcal{V}_A$ and $\mathcal{E}^1_A\bigcup \mathcal{E}^2_A = \mathcal{E}_A$. It follows from $|\mathcal{E}_A| = |\bar{\mathcal{T}}_{\mathcal{L}_n}| \geq 2n-4$ and $|\mathcal{V}_A| = |\mathcal{E}_n| = 2n-3$ that $|\mathcal{E}_A| > |\mathcal{V}_A|-2$. Since the minimum number of edges of a 2-component graph with $n$ vertices is $n-2$, there exists an undirected cycle in $\left(\mathcal{V}_A, \mathcal{E}_A \right)$. Suppose that the undirected cycle is in $\left(\mathcal{V}^1_A, \mathcal{E}^1_A \right)$. Denote it by $\mathcal{C} = (\mathcal{V}_C, \mathcal{E}_C)$, where $\mathcal{V}_C \subseteq \mathcal{V}^1_A$ and $\mathcal{E}_C \subseteq \mathcal{E}^1_A$. Since constraining signed angles $\alpha_{ijk}$ corresponding to $(a_{ji},a_{jk})\in \mathcal{E}_C$ only requires $|\mathcal{V}_C|-1 < |\mathcal{E}_C|$ singed angle constraints, the signed angle constraints given by $\bar{\mathcal{T}}_{\mathcal{L}_n}$ are not minimal, which causes a contradiction. 

            As a result, $\bar{\mathcal{T}}_{\mathcal{L}_n}$ is minimal to satisfy $\mathscr{E}_{\bar{\mathcal{T}}_{\mathcal{L}_n}}(p) = \mathscr{E}_{\mathcal{L}_n}(p)$ only if $|\bar{\mathcal{T}}_{\mathcal{L}_n}| = 2n-4$ and $\bar{\mathcal{T}}_{\mathcal{L}_n}$ is angle connected.
        \end{IEEEproof}
        
        On the one hand, according to the proof of Theorem \ref{lemma: AIS must be connected if and only if the graph is laman graph}, we observe that for a framework $(\mathcal{G},p)$ in $\mathbb{R}^2$, the breadth first search method (BFS) \cite[Section 20.2]{cormen2022introduction} can be applied to extract a spanning tree from the angle index graph $(\mathcal{V}_A(\mathcal{E}), \mathcal{E}_A(\mathcal{T}_{\mathcal{G}}))$. The resulting AIS $\bar{\mathcal{T}}_{\mathcal{G}}$ associated with this spanning tree is angle connected and satisfy $|\bar{\mathcal{T}}_{\mathcal{G}}| = |\mathcal{E}|-1$. On the other hand, according to Theorem \ref{theorem shape fix}-\ref{theorem: ISAR laman subframeworks}, a framework can be uniquely determined by signed angles if and only if it contains an ISAR Laman subframework. Based on these observations, we present Algorithm \ref{alg:1} to construct a minimal GAIS $\bar{\mathcal{T}}_{\mathcal{G}}$ for arbitrary ISAR frameworks, which enables unique shape while satisfying $|\bar{\mathcal{T}}_{\mathcal{G}}| = 2n-4$.

        \begin{algorithm}
		\renewcommand{\algorithmicrequire}{\textbf{Input:}}
		\renewcommand{\algorithmicensure}{\textbf{Output:}}
		\caption{Finding a Minimum GAIS $\bar{\mathcal{T}}_{\mathcal{G}}$ for Any ISAR Framework $(\mathcal{G},p)$}
		\label{alg:1}
		\begin{algorithmic}[1]
			\REQUIRE An ISAR framework $(\mathcal{G}, p)$
			\ENSURE $\bar{\mathcal{T}}_{\mathcal{G}}$
			\STATE Find an ISAR Laman spanning subframework $(\mathcal{L}_n,p)$ from $(\mathcal{G}, p)$ according to Remark \ref{remark :1}. Let $\bar{\mathcal{T}}_{\mathcal{G}} \leftarrow \varnothing$;
			\STATE Find a spanning tree $(\mathcal{V}^{tree}_A,\mathcal{E}^{tree}_A)$ within the angle index graph $(\mathcal{V}_A(\mathcal{L}_n),\mathcal{E}_A(\mathcal{T}_{\mathcal{L}_n}))$ by BFS \cite[Section 20.2]{cormen2022introduction}.
                \FORALL{$(a_{ji},a_{jk})\in\mathcal{E}^{tree}_A$}
                \STATE Let $\bar{\mathcal{T}}_{\mathcal{G}} \leftarrow \{(i,j,k)\}\bigcup\bar{\mathcal{T}}_{\mathcal{G}}$ if $i<k$, $\bar{\mathcal{T}}_{\mathcal{G}} \leftarrow \{(k,j,i)\}\bigcup\bar{\mathcal{T}}_{\mathcal{G}}$ otherwise;
                \ENDFOR
			\STATE \textbf{return} $\bar{\mathcal{T}}_{\mathcal{G}}$
		\end{algorithmic}  
	\end{algorithm}

        \begin{remark}
                Algorithm \ref{alg:1} is a polynomial-time algorithm. As noted in Remark \ref{remark :1}, Line 1 can be achieved via Gaussian elimination  \cite[Section 3.2]{golub2013matrix}, with a time complexity of $\mathcal{O}(|\mathcal{E}|n^2)$ when $m \geq 2n$, and $\mathcal{O}(|\mathcal{E}|^2n)$ otherwise. Since the worst-case time complexity of Lines 2--5 is $\mathcal{O}(n^2)$ (when $(\mathcal{V}_A(\mathcal{L}_n),\mathcal{E}_A(\mathcal{T}_{\mathcal{L}_n}))$ is a complete graph) and $|\mathcal{E}|\leq n(n-1)/2$, the worst-case time complexity of Algorithm \ref{alg:1} is $\mathcal{O}(n^5)$.
        \end{remark}

        \begin{figure}
		\centering
		\includegraphics[width=0.6\linewidth]{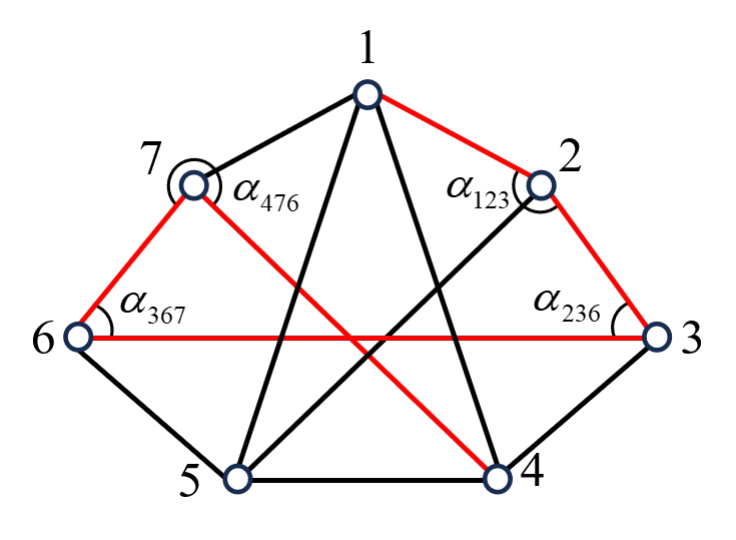}
		\caption{{\footnotesize  Based on the common bearing $b_{12}$, any other bearings can be denoted by the product of a rotation matrix determined by the marked signed angles and $b_{21}$. For example, $b_{47}$ can be depicted as $b_{47} = \mathscr{R}_o(\alpha^*_{47})b_{12}$ where $\alpha^*_{47} = \pi+ \alpha_{123}+\alpha_{236}+\alpha_{367}-\alpha_{476}$.}}
		\label{fig:formation control example}
	\end{figure} 

        \subsection{Minimal RAISs for Signed Angle Rigidity}
        The following result provides an intuitive condition to determine RAISs.
        \begin{lemma} \label{lemma : angle connected set is RAIS}
            For a framework $(\mathcal{G},p)$ in $\mathbb{R}^2$, an AIS $\mathcal{T}^*_{\mathcal{G}}$ is a RAIS if it is angle connected.
        \end{lemma}
        \begin{IEEEproof}
            According to Definition \ref{definition: RAISs}, it is sufficient to show $\rank (R^{\mathcal{T}^*_{\mathcal{G}}}_S(p)) = \rank(R^{\mathcal{T}_{\mathcal{G}}}_S(p))$. Since $\mathcal{T}^*_{\mathcal{G}}$ is angle connected, for any two edges $(i,j), (l,m) \in \mathcal{E}$, one has $b_{lm} = \mathscr{R}_o(\alpha^*_{lm})b_{ij}$ (see Fig. \ref{fig:formation control example} as an example), where 
    \begin{eqnarray}
        \alpha^*_{lm} &\hspace{-2mm}=\hspace{-2mm}& s\pi + \sum_{(i,j,k)\in\mathcal{T}^*_{\mathcal{G}}}c_{ijk}\alpha_{ijk}
    \end{eqnarray}
    for some $c_{ijk},s \in \{\pm1,0\}$. Therefore, the signed angle $\alpha \in [0,2\pi)$ between any pair of bearings in $(\mathcal{G},p)$ satisfies
	\begin{eqnarray} \label{eq: signed angle equation}
		s_1\pi + \alpha &\hspace{-2mm} =\left\{\begin{array}{ll}
			\alpha^* & {\rm if}~ \sin\alpha^* \geqslant 0, \\
			2\pi-\alpha^* & {\rm otherwise},
		\end{array}\right.
	\end{eqnarray}
    where $\alpha^* = s_2\pi + \sum_{(i,j,k)\in\mathcal{T}^*_{\mathcal{G}}}c_{ijk}\alpha_{ijk}$ for some $ s_2, c_{ijk}\in \{0,\pm1\}$ and $s_1 \in \mathbb{N}$. Differentiating both sides of (\ref{eq: signed angle equation}) w.r.t. time leads to 
    \begin{eqnarray}
        \dot{\alpha} &\hspace{-2mm}=\hspace{-2mm}& s_3\sum_{(i,j,k)\in\mathcal{T}^*_{\mathcal{G}}}c_{ijk}\dot{\alpha}_{ijk}
    \end{eqnarray}
    for some $s_3\in\{\pm1\}$, which means any row in $R^{\mathcal{T}_{\mathcal{G}}}_S(p)$ can be represented as a linear combination of the rows indexed by $\mathcal{T}^*_{\mathcal{G}}$. As a result, $\rank (R^{\mathcal{T}^*_{\mathcal{G}}}_S(p)) = \rank(R^{\mathcal{T}_{\mathcal{G}}}_S(p))$.
        \end{IEEEproof}

    Similar to GAISs, the angle connectivity of RAISs is not necessary, see an example in Fig. \ref{fig: non-connected angle index graph}(a), where $\bar{\mathcal{T}}_{\mathcal{G}}$ is a RAIS since $\rank (R^{\bar{\mathcal{T}}_{\mathcal{G}}}_S(p)) = \rank(R^{\mathcal{T}_{\mathcal{G}}}_S(p)) = 1 < 2$ but not angle connected. Nonetheless, we can have the following conclusion.

    \begin{theorem} \label{theorem: AIS is minimal RAIS if and only if the framework is a ISAR laman subframework}
        For an ISAR Laman framework $(\mathcal{L}_n,p)$ in $\mathbb{R}^2$, $\mathcal{T}^*_{\mathcal{L}_n}\subseteq \mathcal{T}_{\mathcal{L}_n}$ is a minimal RAIS if and only if $|\bar{\mathcal{T}}_{\mathcal{L}_n}| = 2n-4$ and $\mathcal{T}^*_{\mathcal{L}_n}$ is angle connected.
    \end{theorem}
    \begin{IEEEproof}
        We first note that for an ISAR framework $(\mathcal{G},p)$, $\mathcal{T}^*_{\mathcal{G}}\subseteq \mathcal{T}_{\mathcal{G}}$ is a RAIS if and only if $\rank (R^{\mathcal{T}^*_{\mathcal{G}}}_S(p)) = \rank(R^{\mathcal{T}_{\mathcal{G}}}_S(p)) = 2n-4$, which implies that $|\mathcal{T}^*_{\mathcal{G}}| \geq 2n-4$.
        
        The sufficiency follows from Lemma \ref{lemma : angle connected set is RAIS} and $|\mathcal{T}^*_{\mathcal{L}_n}| \geq 2n-4$. As for necessity, by employing Lemma \ref{lemma : angle connected set is RAIS} and the techniques used in the proof of Theorem \ref{lemma: AIS must be connected if and only if the graph is laman graph}, it can be verified that the following two cases contradict the fact that $\mathcal{T}^*_{\mathcal{L}_n}$ is minimal to satisfy $\rank (R^{\mathcal{T}^*_{\mathcal{L}_n}}_S(p)) = \rank(R^{\mathcal{T}_{\mathcal{L}_n}}_S(p))$:
        (\rmnum{1}) $|\mathcal{T}^*_{\mathcal{L}_n}| > 2n-4$ and $\mathcal{T}^*_{\mathcal{L}_n}$ is angle connected;
        (\rmnum{2}) $|\mathcal{T}^*_{\mathcal{L}_n}| \geq 2n-4$ and $\mathcal{T}^*_{\mathcal{L}_n}$ is not angle connected.       
    \end{IEEEproof}

        Combining Theorem \ref{lemma: AIS must be connected if and only if the graph is laman graph} with Theorem \ref{theorem: AIS is minimal RAIS if and only if the framework is a ISAR laman subframework}, we obtain the following corollary, which indicates that the AIS generated by Algorithm \ref{alg:1} is both a minimal GAIS and a minimal RAIS.

        \begin{corollary} \label{corollary: summerize minimal RAISs and GAISs}
            For an ISAR Laman framework $(\mathcal{L}_n,p)$ in $\mathbb{R}^2$, the following statements are equivalent:
            \begin{itemize}
                \item[(\rmnum{1})] $\bar{\mathcal{T}}_{\mathcal{L}_n}\subseteq \mathcal{T}_{\mathcal{L}_n}$ is a minimal GAIS;
                \item[(\rmnum{2})] $\bar{\mathcal{T}}_{\mathcal{L}_n}\subseteq \mathcal{T}_{\mathcal{L}_n}$ is a minimal RAIS;
                \item[(\rmnum{3})] $|\bar{\mathcal{T}}_{\mathcal{L}_n}| = 2n-4$ and $\bar{\mathcal{T}}_{\mathcal{L}_n}$ is angle connected.
            \end{itemize}
        \end{corollary}

\begin{figure}
    \centering
    \includegraphics[width=1\linewidth]{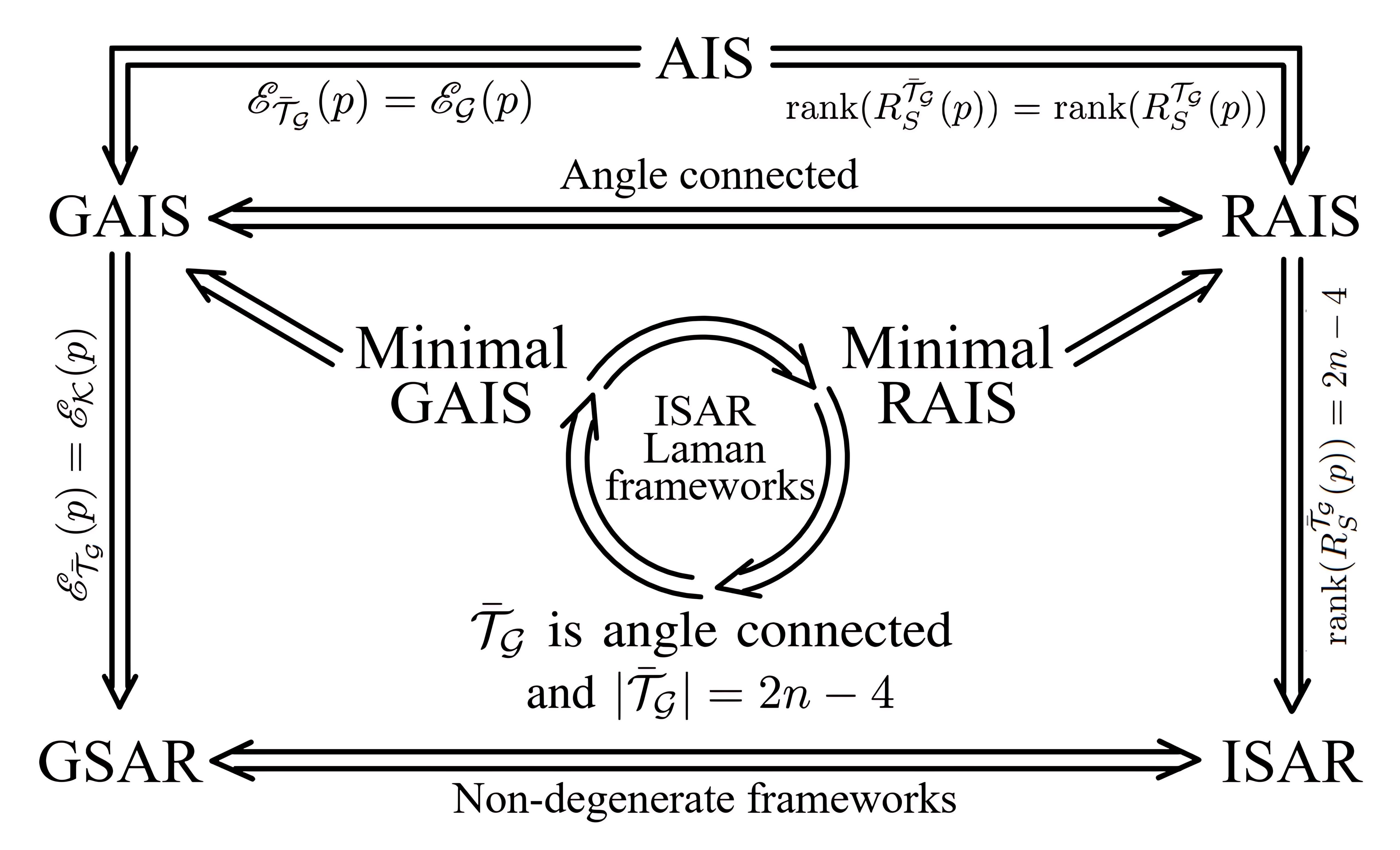} 
    \caption{{\footnotesize A graph illustrates the relationship between different concepts related to AISs and signed angle rigidity properties. An arrow means `implication'. Specifically, `$A\Rightarrow B$ with a condition noted on the arrow' means that $A$ implies $B$ when the condition holds.}}
    \label{AIS and signed angle rigidity}
    \centering
    \end{figure}

    To demonstrate the results stated above, we summarize the relationship between AISs and signed angle rigidity in Fig. \ref{AIS and signed angle rigidity}.

        \section{Signed Angle-Based Sensor Network Localization} \label{sec:5}
        In this section, we study the problem of signed angle-based sensor network localization, where each sensor aims to estimate its own location based on signed angle measurements and communications with neighbors. The results established in the previous sections will provide guidance on how to design localizable networks as well as localization algorithms. 
        
        \subsection{Problem Formulation and Network Localizability}
        Given a sensor network, the sensors can be represented as $\mathcal{V} = \mathcal{A}\bigcup\mathcal{F}$ where the sensors in $\mathcal{A} \subseteq \mathcal{V}$ are called \textit{ anchors} whose locations are known, and the sensors in $\mathcal{F}  = \mathcal{V}\backslash\mathcal{A}$ are called \textit{followers} whose locations are to be determined. Furthermore, the real location of the network is denoted by $p = (p_1,\dots,p_n)^\top$, where $p_i\neq p_j$ for any pair of sensors $i,j\in\mathcal{V}$. The interactions between the sensors are indicated by the fixed and undirected sensing graph $\mathcal{G}=(\mathcal{V},\mathcal{E})$.
        
        We use $(\mathcal{G},p,\mathcal{A})$ to describe a sensor network and consider the SASNL problem as follows.

        \begin{problem} \label{problem 2}
            Consider a sensor network $(\mathcal{G},p,\mathcal{A})$ in $\mathbb{R}^2$, design a distributed algorithm for each follower $i\in\mathcal{F}$ to estimate its location $p_i$ merely using the signed angle measurements $\{\alpha_{ijk}\}_{(i,j,k)\in\mathcal{T}_\mathcal{G}}$ and communicated information from neighbors, i.e., seek a solution to the following nonlinear equations,
    	\begin{eqnarray}
    		\label{eq:problem 2}
    		&\left\{\begin{array}{ll}
    			\frac{\hat{p}_j-\hat{p}_k}{\|\hat{p}_j-\hat{p}_k\|} = \mathscr{R}_o(\alpha_{ijk})\frac{\hat{p}_j-\hat{p}_i}{\|\hat{p}_j-\hat{p}_i\|}  & \forall (i,j,k)\in \mathcal{T}_{\mathcal{G}}, \\
    			\hat{p}_i = p_i & i \in \mathcal{A},
    		\end{array}\right.
    	\end{eqnarray}                 
            where $\hat{p}_i$ is the estimated location of follower $i$.
        \end{problem}

        It should be noted that (\ref{eq:problem 2}) may have multiple solutions, not all of which correspond to the real location of the network. Therefore, it is necessary to investigate the condition for the existence of a unique solution to (\ref{eq:problem 2}), which motivates the introduction of the ``localizability" notion.

        \begin{definition}
            A sensor network $(\mathcal{G},p,\mathcal{A})$ in $\mathbb{R}^2$ is said to be \textit{signed angle localizable} if the real location $p$ is the only feasible solution to (\ref{eq:problem 2}).
        \end{definition}
        
        The following theorem proposes a necessary and sufficient condition for signed angle localizability.
        
        \begin{theorem} \label{theorem: signed angle localizable}
            For an ISAR framework $(\mathcal{G},p)$ in $\mathbb{R}^2$, the sensor network $(\mathcal{G},p,\mathcal{A})$ is signed angle localizable if and only if  $|\mathcal{A}|\geq2 $.
        \end{theorem}
        
        \begin{IEEEproof}
            Sufficiency. Let $q = (q^\top_1,\dots,q^\top_n)^\top\in\mathbb{R}^{2n}$ be a solution to (\ref{eq:problem 2}). Then, we have $q \in \mathscr{E}_{\mathcal{G}}(p)$. It suffices to show $q = p$. Since $(\mathcal{G},p)$ is ISAR, it follows from Theorem \ref{theorem shape fix} that
            \begin{eqnarray}
                q &\hspace{-2mm}=\hspace{-2mm}& c(I_n \otimes\mathscr{R}_o(\theta))p+\textbf{1}_n\otimes\xi,
            \end{eqnarray}
            for some $c\in\mathbb{R}\backslash\{0\}$, $\theta\in[0,2\pi)$, and $\xi\in\mathbb{R}^2$. Without loss of generality, let $1,2$ be the anchors since $|\mathcal{A}|\geq2$. It follows from $q_i = p_i, i\in\mathcal{A}$ that 
            \begin{eqnarray}
                p_1-p_2 &\hspace{-2mm}=\hspace{-2mm}& q_1-q_2 = c\mathscr{R}_o(\theta)(p_1-p_2),
            \end{eqnarray}
            which implies $c\mathscr{R}_o(\theta) = I_2$. Together with the fact that $q_1 = c\mathscr{R}_o(\theta)p_1 + \xi = p_1$, we have $\xi = [0,0]^\top$ and $q = p$. Hence, $(\mathcal{G},p,\mathcal{A})$ is signed-angle-based localizable.

            Necessity. It suffices to show that there exists a solution $q\in\mathbb{R}^{2n}$ of (\ref{eq:problem 2}) such that $q \neq p$ when $|\mathcal{A}|\leq1$. Without loss of generality, suppose $\mathcal{A} = \{1\}$. Let $q = (I_n \otimes\mathscr{R}_o(\frac{\pi}2))p + {\bf 1}_n\otimes \eta$, where $\eta = (I_2-\mathscr{R}_o(\frac{\pi}{2}))p_1$. It can be verified that $q\in\mathscr{E}_{\mathcal{G}}(p)$, $q_1 = p_1$, and $q \neq p$. On the other hand, when $\mathcal{A} = \varnothing$, $q = p + \textbf{1}_n\otimes\xi$ with any vector $\xi\in\mathbb{R}^2$ is a solution of (\ref{eq:problem 2}).
        \end{IEEEproof}
        
\begin{figure}
		\centering
		\includegraphics[width=1\linewidth]{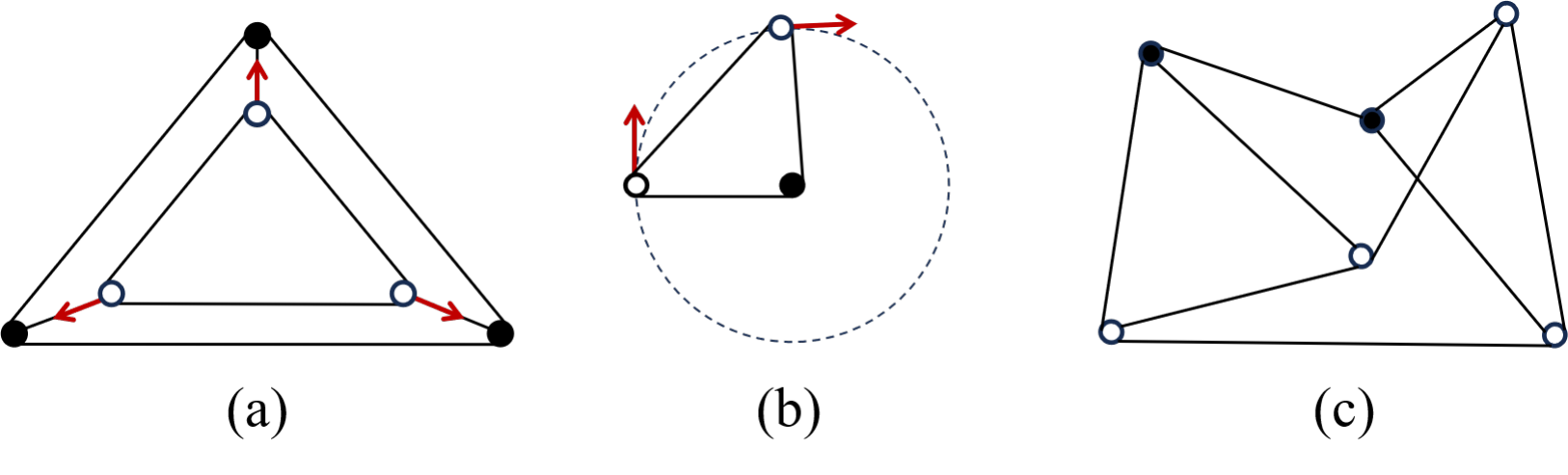}
		\caption{{\footnotesize Examples of sensor networks. Anchors and followers are denoted by solid and hollow dots, respectively. The red arrows represent the infinitesimal signed angle motions corresponding to the followers. (a) and (b) are not signed angle localizable. (c) is signed angle localizable.}}
		\label{fig:localization counterexample}
	\end{figure}

        We provide examples in Fig. \ref{fig:localization counterexample} to show non-localizable networks that violate the condition of Theorem \ref{theorem: signed angle localizable}. Note that the localizability conditions proposed in \cite{chen2022triangular} cannot determine whether the network in Fig. \ref{fig:localization counterexample}(c) is signed angle localizable, as they only apply to trigraphs, which is unnecessary.

        \begin{remark}
            Problem \ref{problem 2} involves all signed angles in the sensor network. However, when $(\mathcal{G},p)$ is ISAR, only partial signed angles are required. That is, $\mathcal{T}_{\mathcal{G}}$ can be replaced by the minimal GAIS $\bar{\mathcal{T}}_{\mathcal{G}}$ induced by Algorithm \ref{alg:1}. When the scenario does not support the centralized implementation of Algorithm \ref{alg:1}, each sensor $i\in\mathcal{V}$ can generate a local AIS ${^i}\bar{\mathcal{T}}_{\mathcal{G}}$ such that the angle index graph $(\mathcal{V}_A(\bar{\mathcal{G}_i}), \mathcal{E}_A({^i}\bar{\mathcal{T}}_{\mathcal{G}}))$ is connected and $|{^i}\bar{\mathcal{T}}_{\mathcal{G}}| = |\mathcal{N}^i_{\mathcal{G}}|-1$, where $\bar{\mathcal{G}}_i = (\bar{\mathcal{V}_i},\bar{\mathcal{E}}_i)$, $\bar{\mathcal{V}_i} = \mathcal{N}^i_{\mathcal{G}}\bigcup\{i\}$, and $\bar{\mathcal{E}}_i = \{(i,j)\in\mathcal{E}:j\in\mathcal{N}^i_\mathcal{G}\}$. Then, the angle index set of the network is $\bar{\mathcal{T}}_{\mathcal{G}} = \bigcup_{i\in\mathcal{V}}{^i}\bar{\mathcal{T}}_{\mathcal{G}}$ with $|\bar{\mathcal{T}}_{\mathcal{G}}| = \sum_{i\in\mathcal{V}}(|\mathcal{N}^i_{\mathcal{G}}|-1) = 2|\mathcal{E}|-n$. It can be verified that $\bar{\mathcal{T}}_{\mathcal{G}}$ is angle connected and sufficient for localization according to Lemma \ref{lemma: angle connected set is GAIS}. However, $\bar{\mathcal{T}}_{\mathcal{G}}$ is not minimal even if $(\mathcal{G},p)$ is minimally ISAR (i.e., $\mathcal{G}$ is a Laman graph) since $|\bar{\mathcal{T}}_{\mathcal{G}}| \geq 2(2n-3)-n =  3n-6 > 2n-4$ when $n>2$.
        \end{remark}   

        \subsection{A Distributed Localization Algorithm}

 Next, we design a distributed algorithm to achieve the localization objective under the following assumption.

        \begin{assumption} \label{ass: 1} 
            The sensor network $(\mathcal{G},p,\mathcal{A})$ considered satisfies: (\rmnum{1}) $(\mathcal{G},p)$ is an ISAR Laman framework; (\rmnum{2}) $|\mathcal{A}| \geq 2$; (\rmnum{3}) there exist two anchors $i,j\in\mathcal{A}$ satisfying $(i,j)\in\mathcal{E}$.
        \end{assumption}
        
        In what follows, we consider that      each sensor $i\in\mathcal{V}$ has access to angle measurements associated with the AIS 
        \begin{eqnarray} \label{eq:angle measurement set}
            {^i}\bar{\mathcal{T}}_{\mathcal{G}} = \{(k,i,j)\in\mathcal{V}^3:(k,i,j)\in\bar{\mathcal{T}}_{\mathcal{G}}\}
        \end{eqnarray}
        where $\bar{\mathcal{T}}_{\mathcal{G}}$ is induced by Algorithm \ref{alg:1} under Assumption \ref{ass: 1}. 
        
        Inspired by \cite{zhao2016localizability}, each follower $i\in\mathcal{F}$ generates a random initial value $\hat{p}_i(0) \in \mathbb{R}^2$ and updates its estimated location by
        \begin{eqnarray} 
            \dot{\hat{p}}_i(t) &\hspace{-2mm}=\hspace{-2mm}&  - \sum_{j\in\mathcal{N}_i}P(\hat{b}_{ij}(t))(\hat{p}_i(t)-\hat{p}_j(t)), \label{eq:location estimators1}
        \end{eqnarray}
        where $\hat{b}_{ij}(t) \in \mathbb{R}^{2}$ is the estimate of the global bearing $b_{ij} \triangleq b_{ij}(p)$ and updated through local communications and signed angle measurements and independent of the position estimates $\hat{p}_i(t)$ for all $i\in\mathcal{V}$ and $t\geq0$, $\hat{p}_j(t)$ is obtained via communications with agent $j\in\mathcal{N}_i$. (In the case of anchors, $\hat{p}_i(t) = p_i$ for all time $t\geq0$.) 

        
        Each sensor $i\in\mathcal{V}$ generates a random initial estimate $\hat{b}_{ij}(0) \in \mathbb{R}^2$ for each $j\in\mathcal{N}_i$ and updates $\hat{b}_{ij}(t)$ by
        \begin{eqnarray}
            \dot{\hat{b}}_{ij}(t) &\hspace{-2mm}=\hspace{-2mm}& - \left[ \sum_{(j,i,k_1)\in\bar{\mathcal{T}}_{\mathcal{G}}}\hspace{-2mm}(\hat{b}_{ij}(t)-\mathscr{R}^\top_o(\alpha_{jik_1})\hat{b}_{ik_1}(t)) \right. \notag \\
            &&\hspace{5mm}+\hspace{-2.5mm}\left. \sum_{(k_2,i,j)\in\bar{\mathcal{T}}_{\mathcal{G}}}\hspace{-2mm}(\hat{b}_{ij}(t)-\mathscr{R}_o(\alpha_{k_2ij})\hat{b}_{ik_2}(t)) \right. \notag \\
            &&\hspace{5mm}+\hspace{-2.5mm}\left. \sum_{(i,j,k_3)\in\bar{\mathcal{T}}_{\mathcal{G}}}\hspace{-2mm}(\hat{b}_{ij}(t)+\mathscr{R}^\top_o(\alpha_{ijk_3})\hat{b}_{jk_3}(t)) \right. \notag \\
            &&\hspace{5mm}+\hspace{-2.5mm}\left. \sum_{(k_4,j,i)\in\bar{\mathcal{T}}_{\mathcal{G}}}\hspace{-2mm}(\hat{b}_{ij}(t)+\mathscr{R}_o(\alpha_{k_4ji})\hat{b}_{jk_4}(t))\right]\hspace{-1mm},\label{eq:location estimators2}
        \end{eqnarray}
        where $\hat{b}_{jk_3}(t)$, $\hat{b}_{jk_4}(t)$, $\alpha_{ijk_3}$, and $\alpha_{k_4ji}$ are obtained via communications with agent $j\in\mathcal{N}_i$. (In the case of $(i,j)\in\mathcal{E}$ and $i,j\in\mathcal{A}$, $\hat{b}_{ij}(t) = b_{ij}$ and $\dot{\hat{b}}_{ij}(t) = 0$ for all $t \geq 0$, where $b_{ij}$ is computed from the locations of anchors $i$ and $j$.)

        Note that the distributed protocol (\ref{eq:location estimators1})--(\ref{eq:location estimators2}) is similar to those in \cite{li2019globally}, but it does not rely on bearing measurements and orientation estimation to get global bearings. Additionally, under the action of (\ref{eq:location estimators2}), $\hat{b}_{ij}(t)$ may not satisfy $\hat{b}_{ij}(t) = \hat{b}_{ji}(t)$ and $\|\hat{b}_{ij}(t)\| = 1$ for all $t\geq0$. Hence, $P(\hat{b}_{ij}(t))$ may not be an orthogonal projection for $t\geq0$. However, these deviations do not affect the convergence of $\hat{p}_i(t)$ to $p_i$ for all $i\in\mathcal{F}$. 

        \begin{theorem}[Global Convergence] \label{theorem:localization theorem}
            Under Assumption \ref{ass: 1} and the AIS (\ref{eq:angle measurement set}), by implementing the algorithm (\ref{eq:location estimators1})--(\ref{eq:location estimators2}),  $\hat{p}_i(t)$ converges to $p_i$ exponentially under arbitrary initial estimates $\hat{p}_i(0)$ and $\{\hat{b}_{ij}(0)\}_{j\in\mathcal{N}_i}$ for all $i\in\mathcal{F}$.
        \end{theorem}
        \begin{IEEEproof}
            See Appendix \ref{proof of Theorem 10}.
        \end{IEEEproof}

\begin{figure}
    \centering
    \includegraphics[width=1\linewidth]{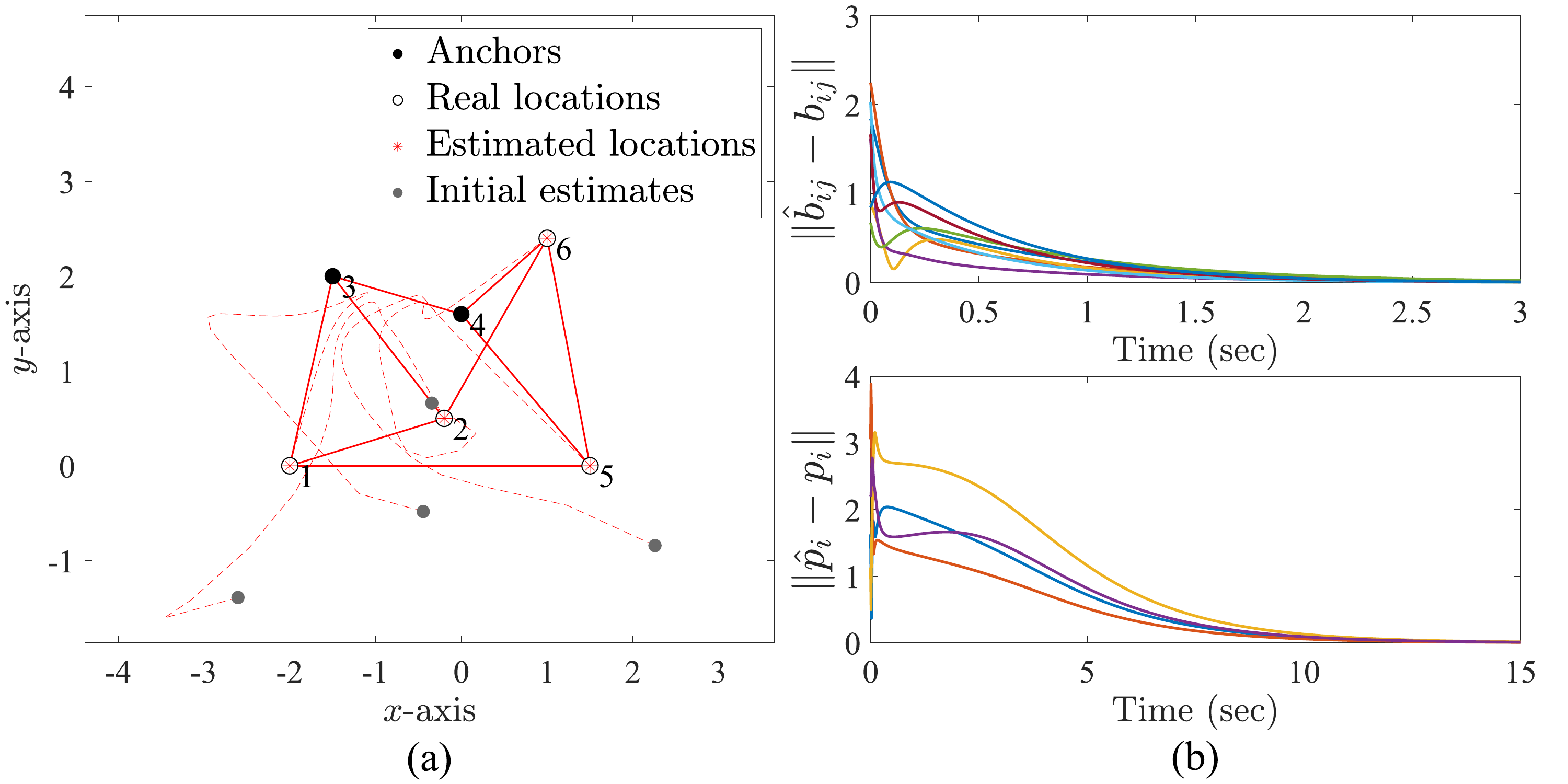} 
    \caption{{\footnotesize (a) Trajectories of the estimation of network locations from randomly initial estates to real positions. (b) Evolution of bearing and location estimation errors.}}
    \centering
    \label{fig:localization example 1}
    \end{figure}

        \begin{example}
        Consider a six-sensor network $(\mathcal{G},p,\mathcal{A})$, which is described in Fig. \ref{fig:localization counterexample}(c) and satisfies Assumption \ref{ass: 1}. The initial estimates $\hat{p}_i(0)$ and $\{\hat{b}_{ij}(0)\}_{j\in\mathcal{N}_i}$ are randomly generated for each sensor $i\in\mathcal{V}$. By implementing (\ref{eq:location estimators1})--(\ref{eq:location estimators2}), the convergence trajectories are shown in Fig. \ref{fig:localization example 1}(a), and the evolutions of bearing and location estimation errors are shown, respectively, in Fig. \ref{fig:localization example 1}(b). From the simulation results, we observe that the estimated locations converge to real locations asymptotically, which is consistent with Theorem \ref{theorem:localization theorem}.
        \end{example}

        \begin{remark}
            The first condition in Assumption \ref{ass: 1} can be relaxed to ISAR frameworks. When the framework $(\mathcal{G},p)$ is ISAR but not a Laman framework, we can first implement Algorithm \ref{alg:1} to induce a minimal GAIS $\bar{\mathcal{T}}_{\mathcal{G}}$. If $\bar{\mathcal{T}}_{\mathcal{G}}$ does not contain a triple $(i,j,k)$ (or $(k,j,i)$) satisfying $(i,j)\in\mathcal{E}$ and $i,j\in\mathcal{A}$, we extend $\bar{\mathcal{T}}_{\mathcal{G}}$ to include such a triple $(i,j,k)$ (or $(k,j,i)$) for some $k\in\mathcal{N}_j$ and $i,j\in\mathcal{A}$. Then, each sensor $i\in\mathcal{V}$ measures the angles indexed by $(\ref{eq:angle measurement set})$, and communicates with neighbors $j\in\bar{\mathcal{N}}_i$ to update the estimates $\hat{b}_{ij}(t)$ and $\hat{p}_i(t)$, where $\bar{\mathcal{N}}_i = \{j,k \in \mathcal{V}:\exists (i,j,k),(k,i,j),$ or $(j,k,i)\in\bar{\mathcal{T}}_{\mathcal{G}}\}$. That is, the subscript ``$j\in\mathcal{N}_i$" in (\ref{eq:location estimators1}) is replaced by ``$j\in\bar{\mathcal{N}}_i$". It can be verified that $\hat{p}_i(t)$ globally converges to $p_i$ for all $i\in\mathcal{F}$ under the distributed protocol (\ref{eq:location estimators1})--(\ref{eq:location estimators2}) when $(\mathcal{G},p)$ is ISAR but not a Laman framework.
        \end{remark}

       \section{Signed Angle Rigid Formation Control} \label{sec:6}

            In this section, we will demonstrate the application of the developed rigidity theory to formation control, where the target formation is described by signed angle constraints and each agent in the formation only measures signed angles to meet these constraints via a distributed control law.
            
            For the group of $n$ agents moving in the plane, we use $\mathcal{V} = \{1,\dots,n\}$ to label the agents. Let $p_i \triangleq p_i(t) \in \mathbb{R}^2$ be the position of agent $i \in \mathcal{V}$, $\beta_i \triangleq \beta_i(t) \in [0,2\pi)$ be the signed angle from the body frame $\sum_i$ to the coordinate frame $\sum_g$, $p \triangleq [p^\top_1,\dots,p^\top_n]$, and $\beta \triangleq [\beta_1,\dots,\beta_n]$. The position and attitude dynamics of agent $i$ is
        \begin{eqnarray}
	       \dot{p}_i &\hspace{-2mm}=\hspace{-2mm}& \mathscr{R}_o(\beta_i)u^p_i(t), \label{equation: position daynamic}\\
              \dot{\beta}_i &\hspace{-2mm}=\hspace{-2mm}& u^a_i(t), \label{equation: attitude daynamic}
        \end{eqnarray}
        where $u^p_i(t)$ and $u^a_i(t)$ are the control input to be designed. 
        
\begin{figure}
		\centering 
		\includegraphics[width=0.7\linewidth]{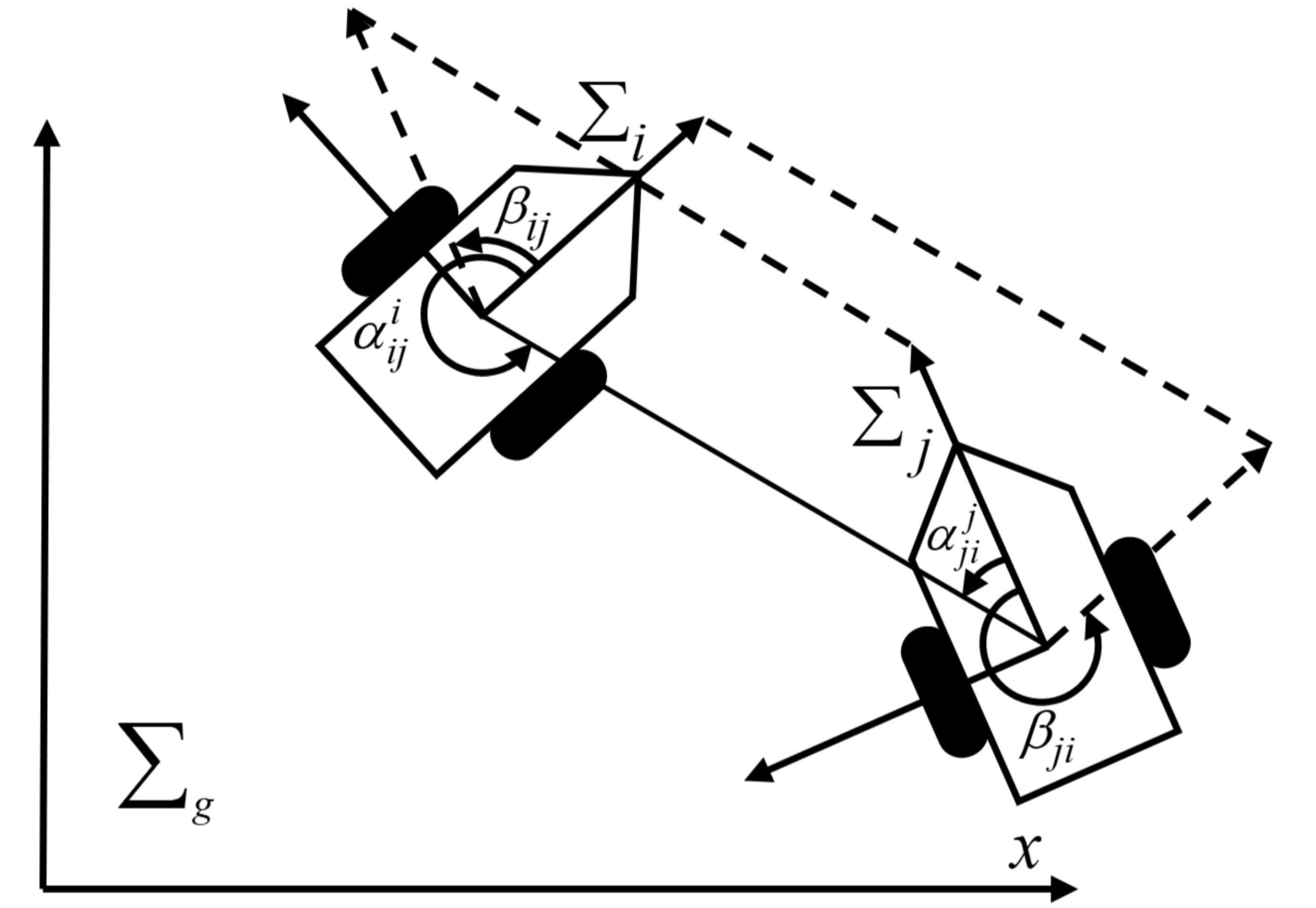}
		\caption{{\footnotesize The illustration of the measured values required for the controllers.}}
		\label{fig:relative orientation}
	\end{figure}

        The sensing relationship between two neighboring agents is represented as an edge in the undirected graph $\mathcal{G}=(\mathcal{V},\mathcal{E})$. That is, agent $i$ can measure the signed angle $\alpha^i_{ij} \in[0,2\pi)$ from the $x$-axis to the edge $(i,j)$ in the body frame $\sum_{i}$, and the signed angles $\beta_{ij} = \beta_{i}-\beta_j$ from $\sum_i$ to $\sum_j$, if $j$ is a neighbor of $i$ in $\mathcal{G}$, see Fig. \ref{fig:relative orientation} for illustrations. We adopt a framework $(\mathcal{G},p^*)$ to denote the target formation. The configuration set that has the same shape as the target formation is described as the set $\mathscr{E}(p^*)$. 
        
        Since ISAR ensures a unique shape of a framework, we consider the following assumption for the target formation.
        
        \begin{assumption}	\label{ass: 2}
	The target formation $(\mathcal{G}, p^*)$ is ISAR.
        \end{assumption}
        
        Note that $\mathscr{E}(p^*)$ is a 4-dimensional smooth manifold under Assumption \ref{ass: 2}. The signed angle-based formation control problem to be solved is stated as below.
\begin{problem}
	Given a set of signed angle constraints $\mathcal{S} = \left\{\alpha_{ijk}(p) = \alpha_{ijk}(p^*), (i,j,k) \in \mathcal{T}_{\mathcal{G}}\right\}$ induced by the framework $(\mathcal{G},p^*)$, design a controller for each agent $i$ based on only the signed angle measurements $\left\{\alpha^i_{ij} \right\}_{j \in \mathcal{N}_i}$ and $\{\beta_{ij}\}_{j \in \mathcal{N}_i}$ such that $\mathscr{E}(p^*)$ is asymptotically stable.
\end{problem}

    Under Assumption \ref{ass: 2}, the graph $\mathcal{G}$ must be connected. Therefore, any bearing in the framework $(\mathcal{G},p^*)$ can be represented as the product of a rotation matrix determined by signed angle constraints and a common bearing, see Fig. \ref{fig:formation control example} for an illustration. As a result, the set $\{\alpha^*_{ij}\}_{j\in\mathcal{N}_i}$ is available for agent $i\in\mathcal{V}$, where $\alpha^*_{ij}$ is the signed angle from the common edge $(1,2)$ to the edge $(i,j)$ in the desired framework $(\mathcal{G},p^*)$, and can be calculated by $\alpha^*_{ij} = s\pi + \sum_{(i,j,k)\in\mathcal{T_{\mathcal{G}}}}c_{ijk}\alpha^*_{ijk}$ for some $s, c_{ijk} \in \{0,\pm1\}$.

        \subsection{The Control Law and Stability Analysis}
        Inspired by the controllers in \cite{zhao2015bearing}, we propose the following signed angle-based control law: 
        \begin{eqnarray}
            u^p_i(t) &\hspace{-2mm}=\hspace{-2mm}& -\sum_{j\in\mathcal{N}_i}\bar{R}(\alpha^i_{ij})\eta(\beta_{ij}, \alpha^*_{ij}),\label{equation:16}\\
            u^a_i(t) &\hspace{-2mm}=\hspace{-2mm}& -\sum_{j\in\mathcal{N}_i}\beta_{ij},\label{equation:17}
        \end{eqnarray}
        where $\eta(\beta_{ij}, \alpha^*_{ij}) = \mathscr{R}_o(\frac{2\alpha^*_{ij}-\beta_{ij}}{2})[0,\cos \frac{\beta_{ij}}{2}]^\top$ and  
        $$\bar{R}(\alpha^i_{ij}) = \begin{bmatrix}
            1-\cos 2\alpha^i_{ij} & -\sin 2\alpha^i_{ij} \\
            -\sin 2\alpha^i_{ij} & 1+\cos 2\alpha^i_{ij}
        \end{bmatrix}.$$ 

        \begin{theorem}[Global Convergence] \label{theorem: global stability}
            Given a group of agents with dynamics (\ref{equation: position daynamic})--(\ref{equation: attitude daynamic}) moving in the plane and the target formation $(\mathcal{G},p^*)$ satisfying Assumption \ref{ass: 2}. By implementing the controller (\ref{equation:16})--(\ref{equation:17}),  $p(t)$ converges to $\mathscr{E}(p^*)$ exponentially under arbitrary initial states $p(0)$ and $\beta(0)$.
    \end{theorem}
    
    \begin{IEEEproof}
    Since $\beta_{ij}$ can be expressed as $\beta_i-\beta_j$, substituting (\ref{equation:17}) into (\ref{equation: attitude daynamic}) leads to
    \begin{eqnarray} 
        \dot{\beta}_i &\hspace{-2mm}=\hspace{-2mm}& \sum_{j\in\mathcal{N}_i}(\beta_j-\beta_i) \label{equation:18}.
    \end{eqnarray}
        It can be observed that (\ref{equation:18}) is the same as the consensus algorithm in \cite[Lemma 1.3]{ren2010distributed} under a connected undirected graph $\mathcal{G}$. Therefore, $\beta_{ij}, (i,j)\in\mathcal{E}$ converges to zero exponentially fast for any initial points $\beta_i(0)\in[0,2\pi), i \in \mathcal{V}$. That is, the body frames of agents will achieve coordinate alignment exponentially. Using similar techniques in the proof of \cite[Theorem 15]{zhao2015bearing}, it can be concluded that the controller (\ref{equation:16}) will drive $b_{ij}(p), (i,j) \in \mathcal{E}$ to either $\mathscr{R}_o(\theta)b_{ij}(p^*)$ (for almost all initial values) or $-\mathscr{R}_o(\theta)b_{ij}(p^*)$ (for the remaining case) exponentially fast, where $\theta \in [0,2\pi)$. It follows that $\alpha_{ijk}(p) \rightarrow \alpha_{ijk}(p^*), (i,j,k) \in \mathcal{T}_{\mathcal{G}}$ as $t \rightarrow \infty$ from any initial values. According to Theorem \ref{theorem shape fix}, the position $p(t)$ of the formation converges to the manifold $\mathscr{E}(p^*)$ globally.
    \end{IEEEproof}

     Theorem \ref{theorem: global stability} implies that the effectiveness of controller (\ref{equation:16})--(\ref{equation:17}) only requires $\mathcal{G}$ to contain a Laman subgraph, which is much milder than the graphical conditions in the literature of angle-based formation control \cite{jing2019angle,chen2020angle}. 
    
\begin{figure}
    \centering
    \includegraphics[width=1\linewidth]{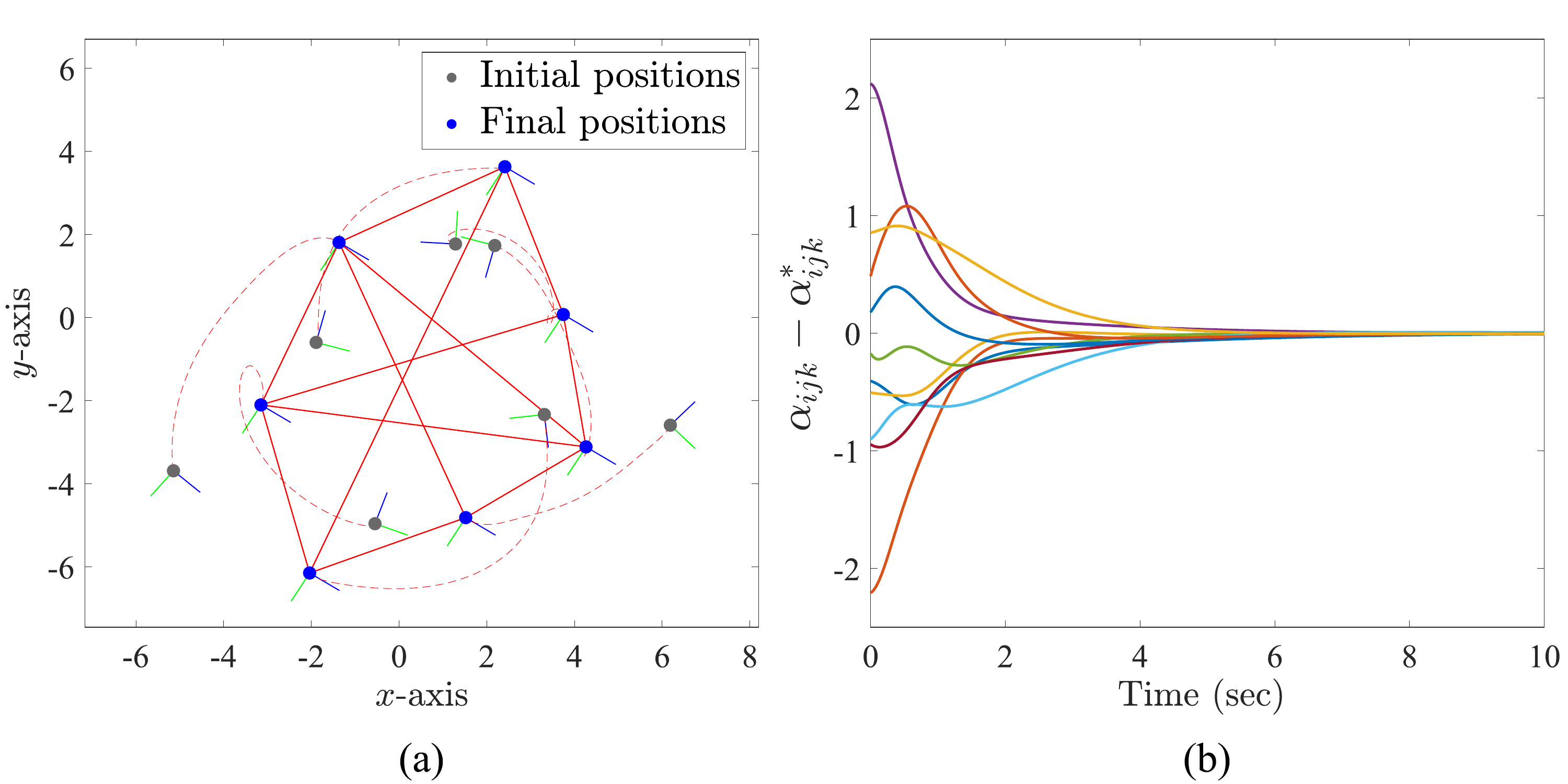} 
    \caption{{\footnotesize (a) Trajectories of seven agents from randomly initial positions to final positions. (b) Evolution of signed angles errors.}}
    \centering
    \label{fig:controller example2}
    \end{figure}

        \begin{example}
        Consider a seven-agent formation, whose target formation shape is described in Fig. \ref{fig:formation control example}. The initial states $\{p_i(0)\}_{i\in\mathcal{V}}$ and ${\beta_i(0)}_{i\in\mathcal{V}}$ are assigned randomly. By implementing the controller (\ref{equation:16})--(\ref{equation:17}), formation trajectories of all agents are shown in Fig. \ref{fig:controller example2}(a), and the evolution of signed angles errors is shown in Fig. \ref{fig:controller example2}(b). From the simulation results, we observe that the formation achieves the desired shape asymptotically, which is consistent with Theorem \ref{theorem: global stability}.
        \end{example}

    \section{Conclusion} \label{sec:7}
        In this paper, signed angle rigidity was reformulated by focusing on graphs. By establishing the complete equivalence between signed angle rigidity and bearing rigidity, we proved that the shape of Laman graphs can be uniquely determined by signed angles for almost all configurations. In order to exclude redundant angle constraints in characterizing signed angle rigidity, an algorithm based on angle index graphs was developed to construct the minimal signed angle constraint set. Building on the developed signed angle rigidity theory, we further designed distributed network localization and formation stabilization algorithms, which are based on signed angle measurements only and achieve global exponential convergence.

        However, the developed signed angle rigidity theory is limited to $\mathbb{R}^2$. Furthermore, the development of this paper relies on establishing the relationship between signed angle rigidity and bearing rigidity. Similar conclusions may hold for angle rigidity in \cite{jing2019angle}, where the angle constraints contain less information. We leave these theoretical extensions as future research directions. 

\section{Appendix: Proofs of Lemmas and Theorems}
\subsection{Proof of Lemma \ref{lemma: new form of signed angle rigidity matrix}} \label{Proof of Lemma lemma: new form of signed angle rigidity matrix}
\textit{Proof: }Differentiating both sides of (\ref{eq:signed angle}) w.r.t. time leads to
        \begin{eqnarray} \label{eq : Lemma1 equation1}
            \hspace{-2mm} \dot{\alpha}_{ijk} \hspace{-2mm}&=&\hspace{-2mm}\frac{b_{jk}^\top P_{ji}}{-\sin\alpha_{ijk}\|e_{ji}\|}\dot{p}_i+\frac{b_{ji}^\top P_{jk}}{-\sin\alpha_{ijk}\|e_{jk}\|}\dot{p}_k \notag \\
            &&-\left(\frac{b_{jk}^\top P_{ji}}{-\sin\alpha_{ijk}\|e_{ji}\|}+\frac{b_{ji}^\top P_{jk}}{-\sin\alpha_{ijk}\|e_{jk}\|}\right)\dot{p}_j.
        \end{eqnarray}
        Note that $e_{ji}$ and $e_{jk}$ can be described by
        \begin{eqnarray}
            e_{ji} = c^x_{ji}b_{jk}+c^y_{ji}b^\bot_{jk}, e_{jk} = c^x_{jk}b_{ji}+c^y_{jk}b^\bot_{ji},
        \end{eqnarray}
        for some $c^x_{ji}, c^y_{ji}, c^x_{jk}, c^y_{jk} \in \mathbb{R}$. Then, one has
        \begin{eqnarray*}
            b^\top_{ji}P_{jk} \hspace{-2mm}&=&\hspace{-2mm} \frac{\left(c^x_{ji}b_{jk}+c^y_{ji}b^\bot_{jk}\right)^\top \hspace{-2mm} P_{jk}}{\|e_{ji}\|}= \sin\alpha_{ijk}(b_{jk})^\top\mathscr{R}_o(\frac{\pi}{2}), \notag \\
            b^\top_{jk}P_{ji} \hspace{-2mm}&=&\hspace{-2mm} \frac{\left(c^x_{jk}b_{ji}+c^y_{jk}b^\bot_{ji}\right)^\top \hspace{-2mm} P_{ji}}{\|e_{jk}\|} = - \sin\alpha_{ijk}(b_{ji})^\top\mathscr{R}_o(\frac{\pi}{2}).
        \end{eqnarray*}
        Together with (\ref{eq : Lemma1 equation1}), we have
        \begin{eqnarray}\label{eq:derivative of signed angle}
            \dot{\alpha}_{ijk} \hspace{-2mm}&=&\hspace{-2mm} \frac{b^\top_{ji}\mathscr{R}_o(\frac{\pi}{2})}{\|e_{ji}\|}\dot{p}_i - \frac{b^\top_{jk}\mathscr{R}_o(\frac{\pi}{2})}{\|e_{jk}\|}\dot{p}_k \notag \\
            &&-\left(\frac{b^\top_{ji}\mathscr{R}_o(\frac{\pi}{2})}{\|e_{ji}\|} - \frac{b^\top_{jk}\mathscr{R}_o(\frac{\pi}{2})}{\|e_{jk}\|}\right)\dot{p}_j.
        \end{eqnarray}
        It can be observed that (\ref{eq:derivative of signed angle}) does not depend on $\sin\alpha_{ijk}$. Hence, the derivative of (\ref{eq:signed angle}) at $\alpha_{ijk} \in [0,2\pi)$ exists and (\ref{signed angle rigid matrix}) can be transformed to (\ref{signed angle rigidity matrix 2}).\hfill{\IEEEQED}

\subsection{Proof of Theorem \ref{corollary ibr isar and iar}} \label{Proof of Theorem 1}

        To prove Theorem \ref{corollary ibr isar and iar}, the following lemma will be used.
        \begin{lemma} \label{lemma:ISAR non-collinear}
            In $\mathbb{R}^2$, a framework $(\mathcal{G},p)$ is ISAR only if for any vertex $j\in\mathcal{V}$, there are at least two distinct neighbors $i,k\in\mathcal{N}_j$ such that $|b^\top_{ji}b_{jk}| \neq 1.$
        \end{lemma}
        \begin{IEEEproof}
            We first show that each vertex has at least two different neighbors. Suppose there is a vertex $j \in \mathcal{V}$ with only one neighbor $i \in \mathcal{V}$.  Let $q = (\textbf{0},\cdots,\textbf{0},(b_{ji})^\top,\textbf{0},\cdots,\textbf{0})^\top \in \mathbb{R}^{2n}$, where $b_{ji}$ consists of the $(2j-1)$-th to $2j$-th components of $q$. According to Lemma \ref{lemma: new form of signed angle rigidity matrix}, for any $\bar{\mathcal{T}}_{\mathcal{G}}\subseteq\mathcal{T}_{\mathcal{G}}$, one has
            \begin{eqnarray}
                R^{\bar{\mathcal{T}}_{\mathcal{G}}}_S(p)q &\hspace{-2mm}=\hspace{-2mm}& \bar{R}_{\bar{\mathcal{T}}_{\mathcal{G}}}(p)\bar{H}q = 0.
            \end{eqnarray}
            Therefore, $q$ is an infinitesimal signed angle motion but not trivial. Hence, a contradiction with the ISAR of $(\mathcal{G},p)$ arises.

            Similarly, if there is a vertex $j$ satisfying $|\mathcal{N}_j| \geq 2$ and $|b^\top_{ji}b_{jk}| = 1$ for any pair of $i,k\in\mathcal{N}_j$, we can find a non-trivial infinitesimal signed angle motion $q = (\textbf{0},\cdots,\textbf{0},(b_{ji})^\top,\textbf{0},\cdots,\textbf{0})^\top \in \mathbb{R}^{2n}$ for some $i \in \mathcal{N}_j$, where $b_{ji}$ consists of the $(2j-1)$-th to $2j$-th components of $q$.
        \end{IEEEproof}

    We next prove Theorem \ref{corollary ibr isar and iar}.

    \textit{Proof: }From the rank condition of ISAR and IBR, it suffices to show that $\rank(R_B(p))$ = $2n-3$ if and only if $\rank(R^{\mathcal{T}_{\mathcal{G}}}_S(p))$ = $2n-4$. 

        According to \cite[Theorem 5]{jing2019angle}, we know that for a framework $(\mathcal{G}, p)$ in $\mathbb{R}^2$, $\rank(R_B(p))=2n-3$ if and only if $\rank(R^{\mathcal{T}_{\mathcal{G}}}_A(p))=2n-4$, where $R^{\mathcal{T}_{\mathcal{G}}}_A(p)$ is the rigidity matrix defined in \cite{jing2019angle} using unsigned angles (i.e., the cosine of the angles) as the entries of the rigidity function\footnote{In this paper, we use the triple $(i,j,k)$ to index the unsigned angle between $b_{ji}$ and $b_{jk}$, which is consistent with the definition in \cite{jing2019angle}.}. Therefore, it suffices to show that $\rank(R^{\mathcal{T}_{\mathcal{G}}}_A(p))=2n-4$ if and only if $\rank(R^{\mathcal{T}_{\mathcal{G}}}_S(p)) = 2n-4$.

        From Lemma \ref{lemma:ISAR non-collinear} and \cite{jing2019angle}, we know that if $(\mathcal{G},p)$ satisfies $\rank(R^{\mathcal{T}_{\mathcal{G}}}_A(p))=2n-4$ or $\rank(R^{\mathcal{T}_{\mathcal{G}}}_S(p)) = 2n-4$, then for any vertex $j\in\mathcal{V}$, there exist at least two distinct neighbors $i,k\in\mathcal{N}_j$ such that $|b^\top_{ji}b_{jk}| \neq 1$. As a result, when the framework satisfies $\rank(R^{\mathcal{T}_{\mathcal{G}}}_A(p))=2n-4$ or $\rank(R^{\mathcal{T}_{\mathcal{G}}}_S(p)) = 2n-4$, we can find an angle connected set $\mathcal{T}^*_{\mathcal{G}} \subseteq \mathcal{T}_{\mathcal{G}}$ by \cite[Algorithm 1]{jing2019angle}, which satisfies $\sin \alpha_{ijk} \neq 0$ for $ (i,j,k)\in \mathcal{T}^*_{\mathcal{G}}$ and $\rank(R^{\mathcal{T}^*_{\mathcal{G}}}_A(p))=\rank(R^{\mathcal{T}_{\mathcal{G}}}_A(p))$. Moreover, with virtue of Lemma \ref{lemma : angle connected set is RAIS} (will be proposed later with its proof independence of this theorem), one has $\rank(R^{\mathcal{T}^*_{\mathcal{G}}}_S(p))=\rank(R^{\mathcal{T}_{\mathcal{G}}}_S(p))$.

        Necessity. We observe from the angle rigidity matrix in \cite{jing2019angle} and equation (\ref{signed angle rigid matrix}) that
		\begin{eqnarray} \label{equation: IAR and ISAR}
			R^{\mathcal{T}^*_{\mathcal{G}}}_S(p) &\hspace{-2mm}=\hspace{-2mm}& {\rm diag}\left(-\frac{1}{{\rm sin}\alpha_{ijk}}\right)R^{\mathcal{T}^*_{\mathcal{G}}}_A(p),
		\end{eqnarray}
        where ${\rm sin}\alpha_{ijk} \neq 0$ for $(i,j,k)\in\mathcal{T}^*_{\mathcal{G}}$. Since $\rank(R^{\mathcal{T}^*_{\mathcal{G}}}_A(p))=\rank(R^{\mathcal{T}_{\mathcal{G}}}_A(p)) =2n-4$, one has 
        \begin{eqnarray}
            \rank(R^{\mathcal{T}^*_{\mathcal{G}}}_S(p)) &\hspace{-2mm}=\hspace{-2mm}& \rank(R^{\mathcal{T}^*_{\mathcal{G}}}_A(p)) = 2n-4.
        \end{eqnarray}
        Together with the fact that $\rank(R^{\mathcal{T}_{\mathcal{G}}}_S(p)) = \rank(R^{\mathcal{T}^*_{\mathcal{G}}}_S(p))$, we have $\rank(R^{\mathcal{T}_{\mathcal{G}}}_S(p)) = 2n-4$.

        Sufficiency. Similar to the discussion in the proof of necessity, it follows from $\rank(R^{\mathcal{T}^*_{\mathcal{G}}}_S(p))=\rank(R^{\mathcal{T}_{\mathcal{G}}}_S(p)) = 2n-4$ and (\ref{equation: IAR and ISAR}) that $\rank(R^{\mathcal{T}_{\mathcal{G}}}_A(p)) = 2n-4$.\hfill{\IEEEQED}
        
\subsection{Proof of Lemma \ref{lemma equivalence between sar and br}} \label{proof of Lemma 3}
        \textit{Proof: }Necessity. Note that for a signed angle $\alpha_{ijk}$, if $\alpha_{ijk}(q) = \alpha_{ijk}(p)$ and $b_{ji}(q) = \mathscr{R}_o(\theta)b_{ji}(p)$ for some $\theta \in [0,2\pi)$, one has $b_{jk}(q) = \mathscr{R}_o(\theta)b_{jk}(p)$. Consider any pair of vertices $i,j \in \mathcal{V}$, there exists an undirected path $\mathcal{P} = (\mathcal{V}_{P},\mathcal{E}_{P})$ between $i$ and $j$ due to the connectivity of $\mathcal{G}$. It follows from $S_{\mathcal{G}}(q) = S_{\mathcal{G}}(p)$ that $b_{l_{k-1}l_k}(q) = \mathscr{R}_o(\theta)b_{l_{k-1}l_k}(p)$ for any $(l_{k-1},l_k) \in \mathcal{E}_P$ and some $\theta \in [0,2\pi)$. Based on the connectivity of $\mathcal{G}$, we have $B_{\mathcal{G}}(q) = \left(I_n \otimes \mathscr{R}_o(\theta) \right) B_{\mathcal{G}}(p)$.

Sufficiency. For any $(i,j),(j,k)\in\mathcal{E}$, it follows from $B_{\mathcal{G}}(q) = \left(I_n \otimes \mathscr{R}_o(\theta) \right) B_{\mathcal{G}}(p)$ that
            \begin{eqnarray} 
			b^{\top}_{jk}(q)b_{ji}(q) \hspace{-2mm}&=&\hspace{-2mm} b^{\top}_{jk}(p)\mathscr{R}^{\top}_o(\theta) \mathscr{R}_o(\theta)b_{ji}(p) \notag\\
            \hspace{-2mm}&=&\hspace{-2mm} b^{\top}_{jk}(p)b_{ji}(p),
		\end{eqnarray}
		\begin{eqnarray} 
			b^{\top}_{jk}(q)\mathscr{R}_o(\frac{\pi}{2})b_{ji}(q) \hspace{-2mm}&=&\hspace{-2mm} b^{\top}_{jk}(p)\mathscr{R}^{\top}_o(\theta)\mathscr{R}_o(\frac{\pi}{2})\mathscr{R}_o(\theta)b_{ji}(p) \notag\\
			&=&\hspace{-2mm} b^{\top}_{jk}(p)\mathscr{R}_o(\frac{\pi}{2})b_{ji}(p).
		\end{eqnarray}
As a result, one has $S_{\mathcal{G}}(q) = S_{\mathcal{G}}(p)$ according to $(\ref{eq:signed angle})$.\hfill{\IEEEQED}

\subsection{Proof of Theorem \ref{theorem shape fix}} \label{proof of Theorem 3}
\textit{Proof: } (\rmnum{1}) $\Leftrightarrow$ (\rmnum{3}). Suppose $(\mathcal{G},p)$ is ISAR. From Lemma \ref{lemma shape fix manifold} and $\mathscr{E}(p)\subseteq\mathscr{E}_{\mathcal{G}}(p)$, it suffices to show that for any $q\in\mathscr{E}_{\mathcal{G}}(p)$, it always holds $q\in\mathscr{E}(p)$. Since $\mathcal{G}$ must be connected, it follows from Lemma \ref{lemma equivalence between sar and br} that $q\in\mathscr{E}_{\mathcal{G}}(p)$ if only if $(I_n \otimes \mathscr{R}_o(\theta))^{-1}q \in B^{-1}_{\mathcal{G}}(B_{\mathcal{G}}(p))$ for some $\theta \in [0,2\pi)$. According to \cite[Theorem 6]{zhao2015bearing}, we know that for a given IBR framework $(\mathcal{G},p)$, one has
        \begin{multline}
            B^{-1}_{\mathcal{G}}(B_{\mathcal{G}}(p)) = \{q\in\mathbb{R}^{2n}:q = cp+\textbf{1}_n \otimes \xi, c\in \mathbb{R}\backslash \{0\}, \\
            \xi \in \mathbb{R}^{2} \}.
        \end{multline}
        Since ISAR is equivalent to IBR, then we have
        \begin{eqnarray}
            (I_n \otimes \mathscr{R}_o(\theta))^{-1}q  &\hspace{-2mm}=\hspace{-2mm}& cp+\textbf{1}_n \otimes \xi,
        \end{eqnarray}
        for some $c\in \mathbb{R}\backslash \{0\}$ and $\xi \in \mathbb{R}^{2}$. As a result, $q\in\mathscr{E}(p)$.
		
		Suppose the shape of $(\mathcal{G},p)$ can be uniquely determined by signed angles up to uniform rotations, translations, and scalings. Then, $p$ must be non-degenerate, and one has $\mathscr{E}(p) = \mathscr{E}_{\mathcal{G}}(p)$ according to Lemma \ref{lemma shape fix manifold}. Therefore, both $\mathscr{E}(p)$ and $\mathscr{E}_{\mathcal{G}}(p)$ are 4-dimentional smooth manifolds. According to \cite[Proposition 3.10]{lee2013introduction}, the tangent space of $\mathscr{E}_{\mathcal{G}}(p)$ at $p$ is 4-dimensional, which indicates that $\rank (R^{\mathcal{T}_{\mathcal{G}}}_S(p)) = 2n-4$. As a result, it follows from Lemma \ref{lemma: trivial motions of signed angle rigidity 2} that $(\mathcal{G},p)$ is ISAR. 	
		
		(\rmnum{2}) $\Leftrightarrow$ (\rmnum{3}). Suppose $(\mathcal{G},p)$ is non-degenerate and GSAR. It suffices to show $\mathscr{E}(p) = \mathscr{E}_{\mathcal{G}}(p)$. Since $\mathscr{E}_{\mathcal{G}}(p) = \mathscr{E}_{\mathcal{K}}(p)$ and (\rmnum{1}) $\Leftrightarrow$ (\rmnum{3}), we just need to show that there exists an ISAR spanning subframework of $(\mathcal{K},p)$. From the non-degeneracy of $(\mathcal{K},p)$, without loss of generality, we can select three non-collinear vertices from $(\mathcal{K},p)$ and index them by $1,2,3$. It can be observed from Lemma \ref{lemma: trivial motions of signed angle rigidity 2} that the framework $(\mathcal{G}^{\prime}, p^\prime)$ embedded by the triangle $\mathcal{G}^{\prime}$ composed of $1,2,3$ is ISAR, where $p^\prime \in \mathbb{R}^{6}$ is the configuration of $1,2,3$ in $p$. Since for any $4\leq j \leq n$, there always exist $i,k\in \{1,2,3\}$ such that $b_{ji}$ and $b_{jk}$ are not collinear. Combining \cite[Theorem 3]{jing2019angle}, Theorem \ref{corollary ibr isar and iar}, and \cite[Lemma 7]{jing2021angle}, we conclude that the framework embedded by the graph constructed by adding vertex $j$ and edges $(j,i)$, $(j,k)$ to $\mathcal{G}^{\prime}$ for $j=4,\cdots n$ iteratively is ISAR.
		
		Suppose $\mathscr{E}(p) = \mathscr{E}_{\mathcal{G}}(p)$. Then, $p$ must be non-degenerate, and hence, it suffices to show that for any $q \in \mathscr{E}_{\mathcal{G}}(p)$, it always holds $q\in\mathscr{E}_{\mathcal{K}}(p)$. Since $\mathscr{E}(p) = \mathscr{E}_{\mathcal{G}}(p)$ and $q \in \mathscr{E}_{\mathcal{G}}(p)$, we have 
        \begin{eqnarray}
            q &\hspace{-2mm}=\hspace{-2mm}& c(I_n \otimes\mathscr{R}_o(\theta))p+\textbf{1}_n \otimes \xi
        \end{eqnarray}
        for some $\theta \in [0,2\pi), c\in \mathbb{R}\backslash \{0\}$, and $\xi \in \mathbb{R}^{2}$. Therefore, one has $B_{\mathcal{K}}(q) = \left(I_n \otimes \mathscr{R}_o(\theta) \right) B_{\mathcal{K}}(p)$ and $S_{\mathcal{K}}(q) = S_{\mathcal{K}}(p)$.\hfill{\IEEEQED}

\subsection{Proof of Theorem \ref{theorem:localization theorem}} \label{proof of Theorem 10}

        To show convergence, similar to \cite[pp.~110]{khalil2002nonlinear}, for the network $(\mathcal{G},p,\mathcal{A})$ and the AIS $\bar{\mathcal{T}}_{\mathcal{G}}$, we define $\hat{\mathcal{G}} \triangleq (\hat{\mathcal{V}},\hat{\mathcal{E}})$ as the bearing graph under (\ref{eq:location estimators2}), which is a directed graph describing the pairwise relationship between bearings in the angle errors feedback according to (\ref{eq:location estimators2}). Here $v_{ij}\in\hat{\mathcal{V}}$ indexes the estimate $\hat{b}_{ij}(t)$ for $i\in\mathcal{V}$ and $j\in\mathcal{N}^i_{\mathcal{G}}$ and is computed by
	\begin{eqnarray}
		\label{eq:signed angle2}
		v_{ij} & =\left\{\begin{array}{ll}
			2 H(i,j,|\mathcal{V}|) & {\rm if}~ i<j, \\
			2 H(i,j,|\mathcal{V}|)+1 & {\rm otherwise},
		\end{array}\right.
	\end{eqnarray} 
        a directed edge $(v_{ik},v_{ij})$ from $v_{ik}$ to $v_{ij}$ exists in $\hat{\mathcal{E}}$ if $\hat{b}_{ik}(t)$ contributes to the update of $\hat{b}_{ij}(t)$, and $(v_{jk},v_{ij})\in \hat{\mathcal{E}}$ if $\hat{b}_{jk}(t)$ is used in the update of $\hat{b}_{ij}(t)$. Let $\bar{\mathcal{A}}$ be the leader set of the angle index graph $(\mathcal{V}_A(\mathcal{G}),\mathcal{E}_A(\bar{\mathcal{T}}_{\mathcal{G}}))$, where $a_{ij}\in\mathcal{V}_A(\mathcal{G})$ belongs to $\bar{\mathcal{A}}$ if $i,j\in\mathcal{A}$. Let $\hat{\mathcal{A}} = \{v_{ij},v_{ji}:a_{ij}\in\bar{\mathcal{A}}\}$ be the leader set of $\hat{\mathcal{G}}$. Under Assumption \ref{ass: 1}, we know that $\hat{\mathcal{A}}\neq \varnothing$.
        
        Similarly to the undirected path, given a directed graph $\hat{\mathcal{G}}$, a directed path from $l_1\in\hat{\mathcal{V}}$ to $l_{n_P}\in\hat{\mathcal{V}}$ is a sequence of directed edges of the form $(l_1,l_2)$,$\dots, (l_{n_P-1},l_{n_P})$. 
        
        Next, we provide a useful lemma with regard to $\hat{\mathcal{G}}$.

        \begin{figure}
		\centering
		\includegraphics[width=0.9\linewidth]{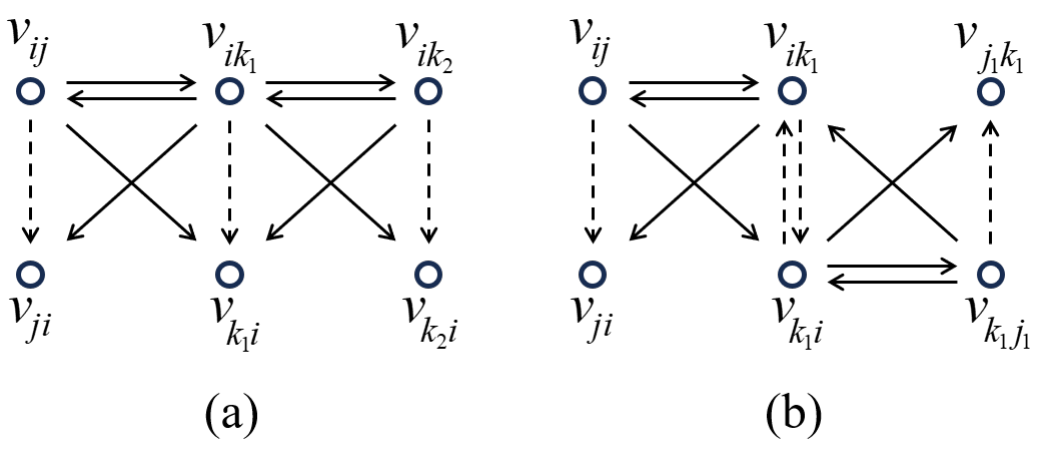}
		\caption{{\footnotesize Illustrations of partial topology relationships of $\hat{\mathcal{G}}$. (a) There exist $(a_{ij},a_{ik_1}),(a_{ik_1},a_{ik_2})\in\mathcal{E}_A$. (b) There exist $(a_{ij},a_{ik_1}),(a_{ik_1},a_{k_1j_1})\in\mathcal{E}_A$. The solid arrow represents the directed edge in $\hat{\mathcal{G}}$ and the dashed arrow represents the existence of a directed path from the starting node to the terminal node.}}
		\label{fig:directed graph with leaders}
	\end{figure}

\begin{lemma} \label{lemma:proof fuzhu of lemma 9}
    Given the network $(\mathcal{G},p,\mathcal{A})$ satisfying Assumption \ref{ass: 1} and the angle connected AIS $\bar{\mathcal{T}}_{\mathcal{G}}$ induced by Algorithm \ref{alg:1}, for any pair of elements $a_{ij},a_{lk}\in\mathcal{V}_A(\mathcal{G})\backslash\bar{\mathcal{A}}$ and the bearing graph $\hat{\mathcal{G}}$ under (\ref{eq:location estimators2}), there exists a directed path from $v_{ij}$ to $v_{lk}$ or from $v_{ji}$ to $v_{lk}$ in $\hat{\mathcal{G}}$ if there exists an undirected path $\mathcal{P}_A = (\mathcal{V}^A_P,\mathcal{E}^A_P)$ satisfying $\mathcal{V}^A_P\bigcap\bar{\mathcal{A}} = \varnothing$ between $a_{ij}$ and $a_{lk}$ in the angle index graph $\mathcal{G}_A(\bar{\mathcal{T}}_{\mathcal{G}}) = \left(\mathcal{V}_A(\mathcal{G}), \mathcal{E}_A(\bar{\mathcal{T}}_{\mathcal{G}}) \right)$.
\end{lemma}
\begin{IEEEproof}
    For simplicity of notation, denote $\mathcal{G}_A(\bar{\mathcal{T}}_{\mathcal{G}}) = \left(\mathcal{V}_A(\mathcal{G}), \mathcal{E}_A(\bar{\mathcal{T}}_{\mathcal{G}}) \right)$ by $\mathcal{G}_A=\left(\mathcal{V}_A, \mathcal{E}_A \right)$. We first note that for any pair of $a_{ij},a_{ik_1}\in\mathcal{V}_A\backslash\bar{\mathcal{A}}$ with $(a_{ij},a_{ik_1})\in\mathcal{E}_A$ , it follows from (\ref{eq:location estimators2}) and the definition of $\hat{\mathcal{E}}$ that the directed edges $(v_{ij},v_{ik_1})$, $(v_{ik_1},v_{ij})$, $(v_{ik_1},v_{ji})$, and $(v_{ij},v_{k_1i})$ belong to $\hat{\mathcal{E}}$ (see Fig. \ref{fig:directed graph with leaders} for examples). It also implies that $v_{ij}$ and $v_{ik_1}$ are mutually reachable in $\hat{\mathcal{G}}$ and there are two directed paths from $v_{ij}$ to $v_{ji}$ and from $v_{ik_1}$ to $v_{k_1i}$, respectively. 
    
    Given any pair of $a_{ij},a_{lk}\in\mathcal{V}_A\backslash\bar{\mathcal{A}}$. If $a_{ij}$ is adjacent to $a_{lk}$ in $\mathcal{G}_A$, then the conclusion holds, immediately. If $a_{ij}$ is not adjacent to $a_{lk}$ in $\mathcal{G}_A$, then there exists an undirected path $\mathcal{P}_A = (\mathcal{V}^A_P,\mathcal{E}^A_P)$ between $a_{ij}$ and $a_{lk}$, where $\mathcal{V}^A_P\subseteq\mathcal{V}_A\backslash\bar{\mathcal{A}}$, $\mathcal{E}^A_P\subseteq\mathcal{E}_A$, and $|\mathcal{V}^A_P|\geq 3$. Without loss of generality, suppose that $(a_{ij},a_{ik_1})\in\mathcal{E}^A_P$. The rest of the proof is to show that if $a_{gh}$ is a neighbor of $a_{ik_1}$ (other than $a_{ij}$) in $\mathcal{P}_A$, there exists two directed paths from $v_{ij}$ to $v_{gh}$ and to $v_{hg}$, respectively. Then, the statement can be proved by induction.
    
    Since the number of neighbors of $a_{ik_1}$ in $\mathcal{P}_A$ is two, one of the following two cases holds:

    Case 1: $\exists (a_{ik_1},a_{ik_2}) \in \mathcal{E}^A_P$. Since $(a_{ij},a_{ik_1}),(a_{ik_1},a_{ik_2})\in\mathcal{E}_A$, it follows that for any pair of vertices in $\{v_{ij},v_{ik_1},v_{ik_2}\}$, they are mutually reachable in $\hat{\mathcal{G}}$ (see Fig. \ref{fig:directed graph with leaders}(a)), and there are two directed paths from $v_{ij}$ to $v_{ik_2}$ and to $v_{k_{2}i}$, respectively. 

    Case 2: $\exists (a_{ik_1},a_{k_1j_1}) \in \mathcal{E}^A_P$. It follows from $(a_{ij},a_{ik_1}),(a_{ik_1},a_{k_1j_1})\in\mathcal{E}_A$ that for any pair of vertices in $\{v_{ik_1},v_{k_1i},v_{ij},v_{k_1j_1}\}$, they are mutually reachable in $\hat{\mathcal{G}}$ (see Fig. \ref{fig:directed graph with leaders}(b)), and there are two directed paths from $v_{ij}$ to $v_{k_1j_1}$ and to $v_{j_1k_1}$, respectively.
   \end{IEEEproof}

    To prove Theorem \ref{theorem:localization theorem}, the following lemma will be used.
        
        \begin{lemma} \label{lemma : bijbji}
            Given the network $(\mathcal{G},p,\mathcal{A})$ satisfying Assumption \ref{ass: 1} and the angle connected AIS $\bar{\mathcal{T}}_{\mathcal{G}}$ induced by Algorithm \ref{alg:1}, there exists a path from some vertex in $\hat{\mathcal{A}}$ to $v_{lk}$ for any vertex $v_{lk}\in \hat{\mathcal{V}}\backslash\hat{\mathcal{A}}$ under (\ref{eq:location estimators2}).
        \end{lemma}
        \begin{IEEEproof}
            Let the angle index graph of $\mathcal{G}$ and $\bar{\mathcal{T}}_{\mathcal{G}}$ be $\mathcal{G}_A = (\mathcal{V}_A,\mathcal{E}_A)$. The rest of the proof is to show that under the cases (\rmnum{1}) $|\bar{\mathcal{A}}| = 1$ and (\rmnum{2}) $|\bar{\mathcal{A}}|\geq2$, the statement holds.
            
            When $|\bar{\mathcal{A}}| = 1$, without loss of generality, let $\bar{\mathcal{A}} = \{a_{12}\}$, and suppose $(a_{12},a_{1k})\in\mathcal{E}_A$. Since $(\mathcal{V}_A,\mathcal{E}_A)$ is a tree, we consider the following three cases:

            Case 1: $a_{12}$ is a leaf vertex and $\exists(a_{1k},a_{1k_2})\in\mathcal{E}_A$. Since $(v_{12},v_{1k})\in\hat{\mathcal{E}}$, it suffices to show that for any $v_{ij}\in \hat{\mathcal{V}}\backslash\hat{\mathcal{A}}$, there exists a directed path from $v_{1k}$ to $v_{ij}$. Recall that $a_{1k},a_{1k_2}\notin\bar{\mathcal{A}}$ and the graph $\bar{\mathcal{G}}_A = (\bar{\mathcal{V}}_A,\bar{\mathcal{E}}_A)$ formed by removing $a_{12}$ from $\mathcal{G}_A$ is still a connected graph. Hence, there exists an undirected path $\mathcal{P}_A = (\mathcal{V}^A_P,\mathcal{E}^A_P)$ between $a_{1k}$ and $a_{ij}$ in the angle index graph $\mathcal{G}_A$, where $\mathcal{V}^A_P\subseteq\mathcal{V}_A\backslash\bar{\mathcal{A}}$ and $\mathcal{E}^A_P\subseteq\mathcal{E}_A$. It follows from Lemma \ref{lemma:proof fuzhu of lemma 9} that there is a directed path from $v_{1k}$ to $v_{ij}$ or from $v_{k1}$ to $v_{ij}$ in $\hat{\mathcal{G}}$. Since $(a_{1k},a_{1k_2}) \in\mathcal{E}_A$ and $a_{1k},a_{1k_2}\notin\bar{\mathcal{A}}$, there is a directed path from $v_{1k}$ to $v_{k1}$ and a directed path from $v_{1k}$ to $v_{ij}$.

            Case 2: $a_{12}$ is a leaf vertex and $\exists(a_{1k},a_{kj_1})\in\mathcal{E}_A$. Using techniques similar to Case 1, it can be verified that for any $v_{ij}\in \hat{\mathcal{V}}\backslash\hat{\mathcal{A}}$, there is a directed path from $v_{12}$ to $v_{ij}$.

            Case 3: $a_{12}$ is a branch vertex. Let $\mathcal{N}_{a_{12}} = \{a^1_{12},\dots,a^m_{12}\}$ be the neighbor set of $a_{12}$ in $\mathcal{G}_A$. Since $a_{12}$ is a branch vertex, $|\mathcal{N}_{a_{12}}| = m \geq 2$ and we can find $m$ subtrees $\bar{\mathcal{G}}^i_A = \left(\bar{\mathcal{V}}^i_A, \bar{\mathcal{E}}^i_A \right)$ (for $i = 1,\dots,m$) from the tree $\mathcal{G}_A$, where $\bar{\mathcal{V}}^i_A = \mathcal{V}^i_A\bigcup\{a_{12}\}$, $\mathcal{V}^i_A$ is the set of the vertices $a_{uv}\in\mathcal{V}_A$ that have a path from $a_{12}$ to $a_{uv}$ passing through $a^i_{12}$, and $\mathcal{E}^i_A = \{(a_{uv},a_{ef})\in\mathcal{E}_A:a_{uv},a_{ef}\in\mathcal{V}^i_A\}$. Note that $a_{12}$ is a leaf vertex in each subtree $\mathcal{G}^i_A$. Therefore, it follows from Case 1--2 that for any $a_{jk}\in \bar{\mathcal{V}}^i_A\backslash\{a_{12}\}$, there exist two directed paths from $v_{12}$ to $v_{jk}$ and to $v_{kj}$, respectively.

            When $|\bar{\mathcal{A}}|\geq2$, one can form the graph $\bar{\mathcal{G}}_A = \left(\bar{\mathcal{V}}_A, \bar{\mathcal{E}}_A \right)$ by removing the vertices $a_{ij}\in\bar{\mathcal{A}}$ from $\mathcal{G}_A$. Let $m$ be the number of connected components in $\bar{\mathcal{G}}_A$, and denote each connected component by $\bar{\mathcal{G}}^i_A = \left(\bar{\mathcal{V}}^i_A, \bar{\mathcal{E}}^i_A \right)$ (for $i = 1,\dots,m$). Since each component $\bar{\mathcal{G}}^i_A$ contains a vertex $a_{lk}\in\bar{\mathcal{V}}^i_A$ that is adjacent to some vertex $a_{uv}\in\bar{\mathcal{A}}$ in $\mathcal{G}_A$, we can further make $\bar{\mathcal{V}}^i_A = \bar{\mathcal{V}}^i_A\bigcup\{a_{uv}\}$ and $\bar{\mathcal{E}}^i_A = \bar{\mathcal{E}}^i_A\bigcup\{(a_{uv},a_{lk})\}$. Then $\bar{\mathcal{G}}^i_A$ becomes a tree after this augmentation and has only one vertex $a_{uv}\in\bar{\mathcal{V}}^i_A$ belonging to $\bar{\mathcal{A}}$. The analysis of the case $|\bar{\mathcal{A}}| = 1$ shows that for any vertex $a_{lk}\in \bar{\mathcal{V}}^i_A\backslash\bar{\mathcal{A}}$, there exist two directed path from some vertex in $\hat{\mathcal{A}}$ to $v_{lk}$ and to $v_{kl}$, respectively. 
        \end{IEEEproof}

         We next prove Theorem \ref{theorem:localization theorem}.
        \begin{IEEEproof}
            Note that according to Theorem \ref{corollary ibr isar and iar} and \cite[Corollary 3]{zhao2016localizability}, Assumption \ref{ass: 1} implies that $(\mathcal{G},p,\mathcal{A})$ is bearing localizable \cite{zhao2016localizability}. Therefore, we first prove that $\hat{b}_{ij}(t)$ exponentially converges to the global bearing $b_{ij}$ for all $j \in \mathcal{N}_i$. Then, the proof of the convergence of $\hat{p}_i(t)$ to $p_i$ for all $i\in\mathcal{F}$ is similar to that in \cite[Theorem 2]{li2019globally}.

            Let $B_{ij} \triangleq [b_{ij},b^\bot_{ij}]$. Given any $\alpha_{ijk} \in [0,2\pi)$, it follows from $\mathscr{R}_o(\alpha_{jik})b_{ij} = b_{ik}$ and $\mathscr{R}_o(\alpha_{jik})b^\bot_{ij} = b^\bot_{ik}$ that $\mathscr{R}_o(\alpha_{jik})B_{ij} = B_{ik}$. Therefore, multiplying $B^\top_{ij}$ on both sides of (\ref{eq:location estimators2}), one has
            \begin{eqnarray} 
            B^\top_{ij}\dot{\hat{b}}_{ij} &\hspace{-2mm}=\hspace{-2mm}& - \left[ \sum_{(j,i,k_1)\in\bar{\mathcal{T}}_{\mathcal{G}}}(B^\top_{ij}\hat{b}_{ij}-B^\top_{ik_1}\hat{b}_{ik_1}) \right. \notag \\
            &&\hspace{5mm}+\hspace{-2.5mm}\left. \sum_{(k_2,i,j)\in\bar{\mathcal{T}}_{\mathcal{G}}}(B^\top_{ij}\hat{b}_{ij}-B^\top_{ik_2}\hat{b}_{ik_2}) \right. \notag \\
            &&\hspace{5mm}+\hspace{-2.5mm}\left. \sum_{(i,j,k_3)\in\bar{\mathcal{T}}_{\mathcal{G}}}(B^\top_{ij}\hat{b}_{ij}-B^\top_{jk_3}\hat{b}_{jk_3}) \right. \notag \\
            &&\hspace{5mm}+\hspace{-2.5mm}\left. \sum_{(k_4,j,i)\in\bar{\mathcal{T}}_{\mathcal{G}}}(B^\top_{ij}\hat{b}_{ij}-B^\top_{jk_4}\hat{b}_{jk_4})\right], \label{eq:location estimators7}
        \end{eqnarray}
        where we use the fact that $B_{ij} = -B_{ji}$. Let $x_{ij} \triangleq B^\top_{ij}\hat{b}_{ij}$ for $j \in \mathcal{N}_i$ and $\mathcal{N}_{v_{ij}}$ be the neighbor set of the vertex $v_{ij}$ in the bearing graph $\hat{\mathcal{G}}$ under (\ref{eq:location estimators2}). Then, (\ref{eq:location estimators7}) can be rewritten as 
        \begin{eqnarray}\label{eq:consensus algorithm with multiple leaders}
            \dot{x}_{ij} &\hspace{-2mm}=\hspace{-2mm}& -\left[\sum_{v_{ik}\in\mathcal{N}_{v_{ij}}}\hspace{-2mm}(x_{ij}-x_{ik})+\hspace{-2mm}\sum_{v_{jk}\in\mathcal{N}_{v_{ij}}}\hspace{-2mm}(x_{ij}-x_{jk})\right].
        \end{eqnarray}
        Since for each vertex $v_{ij}\in \hat{\mathcal{V}}\backslash\hat{\mathcal{A}}$, there exists a directed path from some vertex in $\hat{\mathcal{A}}$ to $v_{ij}$ by Lemma \ref{lemma : bijbji}, and $x_{jk} = B^\top_{jk}\hat{b}_{jk} = B^\top_{jk}b_{jk} = [1,0]^\top$ for $v_{jk}\in \hat{\mathcal{A}}$ and all $t\geq0$, $x_{ij}$ converges to $[1,0]^\top$ for all $j \in \mathcal{N}_i$ according to the result in \cite[Theorem 5.1]{ren2010distributed}. It implies that the system converges to
            \begin{equation}
                \left\{
                \begin{array}{l}
                    \hat{b}^\top_{ij}b_{ij} =  \|\hat{b}_{ij}\|\cos\theta = 1, \\
                    \hat{b}^\top_{ij}b^\bot_{ij} =  \|\hat{b}_{ij}\|\sin\theta = 0,
                \end{array}
                \right.
                \label{eq:myequation}
            \end{equation}
        where $\theta$ is the signed angle from $\hat{b}_{ij}$ to $b_{ij}$. It follows from (\ref{eq:myequation}) that $\|\hat{b}_{ij}\|\mathscr{R}_o(\theta)[1,0]^\top = [1,0]^\top$. Therefore, one has $\|\hat{b}_{ij}\| = 1$ and $\theta = 0$. As a result, $\hat{b}_{ij}\rightarrow b_{ij}$ as $t \rightarrow \infty$ for all $j \in \mathcal{N}_i$. Since the system (\ref{eq:consensus algorithm with multiple leaders}) is linear, the convergence is exponentially fast \cite[Corollary 4.3]{khalil2002nonlinear}.
        \end{IEEEproof}

\bibliographystyle{unsrt}  
{\footnotesize
\bibliography{references.bib}  
}

\end{document}